\documentclass[reqno]{amsart}

\usepackage{graphicx}
\usepackage{amsfonts,amsmath, amssymb}

\numberwithin{equation}{section}

\usepackage{hyperref}
\usepackage{marginnote}

\newcommand{\dis}{\displaystyle}

\usepackage{color}

\newcommand{\R}{\mathbb{R}}

\newtheorem{theorem}{Theorem}[section]
\newtheorem{corollary}[theorem]{Corollary}
\newtheorem{lemma}[theorem]{Lemma}
\newtheorem{proposition}[theorem]{Proposition}
\newtheorem{remark}[theorem]{Remark}

\def\v{\varepsilon}

\def\t{\theta}

\def\a{\alpha}
\def\b{\beta}
\def\g{\gamma}
\def\d{\delta}
\def\l{\lambda}

\def\r{\rho}
\def\s{\sigma}
\def\z{\zeta}

\def\f{\frac}

\begin{document}

\title[The Boltzmann equation with large amplitude initial data]{
Global well-posedness of the Boltzmann equation with large amplitude initial data}

\author[R.-J. Duan]{Renjun Duan}
\address[R.-J. Duan]{Department of Mathematics, The Chinese University of Hong Kong, Hong Kong}
\email{rjduan@math.cuhk.edu.hk}

\author[F.-M. Huang]{Feimin Huang}
\address[F.-M. Huang]{Institute of Applied Mathematics, AMSS, CAS, Beijing 100190, P.R.~China}
\email{fhuang@amt.ac.cn}

\author[Y. Wang]{Yong Wang}
\address[Y. Wang]{Institute of Applied Mathematics, AMSS, CAS, Beijing 100190, P.R.~China}
\email{yongwang@amss.ac.cn}

\author[T. Yang]{Tong Yang}
\address[T. Yang]{Department of Mathematics, City University of Hong Kong, Hong Kong}
\email{matyang@cityu.edu.hk}

\date{\today}

\maketitle

\begin{abstract}
The global well-posedness of the Boltzmann equation with initial data of large amplitude has remained a long-standing open problem. In this paper, by developing a new $L^\infty_xL^1_{v}\cap L^\infty_{x,v}$  approach,   we prove the global existence and uniqueness of mild solutions to the Boltzmann equation in the whole space or torus for a class of initial data with bounded velocity-weighted $L^\infty$ norm under some  smallness condition on $L^1_xL^\infty_v$ norm as well as  defect mass, energy and entropy so that the initial data  allow  large amplitude oscillations. Both the hard and  soft potentials with angular cut-off are  considered,  and the large time behavior of solutions in $L^\infty_{x,v}$ norm with explicit rates of convergence is also studied. 

\end{abstract}

\tableofcontents

\thispagestyle{empty}

\section{Introduction}

In this paper, we consider the Boltzmann equation
\begin{equation}\label{1.1}
	F_t+v\cdot\nabla_x F=Q(F,F),
\end{equation}
where $F(t,x,v)\geq 0$ is the density distribution function for the gas particles with position $x\in\Omega=\mathbb{R}^3$ or $\mathbb{T}^3$ and velocity $v\in\mathbb{R}^3$ at  time $t>0$. The Boltzmann collision term $Q(F,F)$ on the right is defined in terms of the following bilinear form 
\begin{multline}\label{1.2}
	Q(F_1,F_2)\equiv\int_{\mathbb{R}^3}\int_{\mathbb{S}^2} B(v-u,\t)F_1(u')F_2(v')\,{d\omega du}\\
	-\int_{\mathbb{R}^3}\int_{\mathbb{S}^2} B(v-u,\t)F_1(u)F_2(v)\,{d\omega du}:=Q_+(F_1,F_2)-Q_-(F_1,F_2),
\end{multline}
where the relationship between the post-collison velocity $(v',u')$ of two particles with the pre-collision velocity $(v,u)$ is given by
\begin{equation*}
	u'=u+[(v-u)\cdot\omega]\omega,\quad v'=v-[(v-u)\cdot\omega]\omega,
\end{equation*}
for $\omega\in \mathbb{S}^2$, which can be determined by conservation laws of momentum and energy: 
\begin{equation*}
u'+v'=u+v,\quad |u'|^2+|v'|^2=|u|^2+|v|^2.
\end{equation*}
The Boltzmann collision kernel $B=B(v-u,\theta)$ in \eqref{1.2} depends only on $|v-u|$ and $\theta$ with $\cos\theta=(v-u)\cdot \omega/|v-u|$.   Throughout this paper,  we  consider both the hard and soft potentials under the Grad's angular cut-off assumption, for instance,
\begin{equation}\label{1.4}
	B(v-u,\t)=|v-u|^{\gamma}b(\t),
\end{equation}
with 
\begin{equation}
-3<\gamma\leq 1,\quad 0\leq b(\t)\leq C|\cos\t|,\notag
\end{equation}
for a postive constant $C>0$. We consider the Boltzmann equation \eqref{1.1} with the following  initial data
\begin{equation}\label{1.5-1}
	F(t,x,v)|_{t=0}=F_0(x,v).
\end{equation}

To look for a solution $F(t,x,v)$ to the Cauchy problem \eqref{1.1} and \eqref{1.5-1}, let us take a reference global Maxwellian
\begin{equation}
\mu(v)=\f{1}{(2\pi)^{\f32}}\exp\left(-\f{|v|^2}{2}\right),\notag
\end{equation}
which is normalized to have unit density, zero bulk velocity and unit temperature. Formally, as introduced in \cite{Guo},  $F(t,x,v)$ 
satisfies the conservations laws of defect mass, momentum, energy:
\begin{align}
\int_{\Omega}\int_{\mathbb{R}^3}(F(t,x,v)-\mu(v))dvdx&= \int_{\Omega}\int_{\mathbb{R}^3}(F_0(x,v)-\mu(v))dvdx:=
M_0,\label{1.13}\\
\int_{\Omega}\int_{\mathbb{R}^3}v(F(t,x,v)-\mu(v))dvdx&= \int_{\Omega}\int_{\mathbb{R}^3}v(F_0(x,v)-\mu(v))dvdx:=
J_0,\label{1.13-1}\\
\int_{\Omega}\int_{\mathbb{R}^3}|v|^2(F(t,x,v)-\mu(v))dvdx&= \int_{\Omega}\int_{\mathbb{R}^3}|v|^2(F_0(x,v)-\mu(v))dvdx:=
E_0,\label{1.13-2}
\end{align}
as well as the inequality of defect entropy 
\begin{equation}\label{1.14-1}
\int_{\Omega}\int_{\mathbb{R}^3}\Big\{F(t,x,v)\ln{F(t,x,v)}-\mu(v)\ln{\mu(v)}\Big\}dvdx
\leq \int_{\Omega}\int_{\mathbb{R}^3}\Big\{F_0\ln{F_0}-\mu(v)\ln{\mu(v)}\Big\}dvdx.
\end{equation}
By defining
\begin{equation}
	\notag
	\mathcal{E}(F(t)):=
	\int_{\Omega}\int_{\mathbb{R}^3}\Big\{F(t,x,v)\ln F(t,x,v)-\mu\ln\mu\Big\} dvdx+\left[\f{3}{2}\ln(2\pi)-1\right]M_0+\f12E_0,
\end{equation}
it follows by a direct calculation 
that
  \begin{equation}
  \notag
\mathcal{E}(F(t))
\geq 0,
  \end{equation}
for all $t\geq 0$. Note, in particular, that $\mathcal{E}(F_0)\geq 0$ holds true for any function $F_0(x,v)\geq 0$.

\medskip

The Boltzmann equation is a fundamental model in the collisional kinetic theory, and there are  enormous literatures on its well-posedness theories, cf.~\cite{CIP} and \cite{Vi} and the references therein. Among them we mention some works only in the spatially inhomogeneous framework; for the spatially homogenous Boltzmann equation, interested readers may refer to Carleman \cite{Car}  as well as the recent work \cite{LM15} and references therein.  For general initial data in $L^\infty$ framework, the local existence and uniqueness was firstly investigated by Kaniel and Shinbrot \cite{Kaniel-S} and the global existence was later obtained by Illner and Shinbrot \cite{IIIner-S} under additional smallness assumption on velocity weighted $L^\infty$ norm. 
It is well known that for general initial data with finite mass, energy and entropy, the global existence of renormalized solutions was proved  by DiPerna and Lions \cite{D-Lion}; 
the uniqueness of such solutions, however,  is unknown. Moreover, the convergence of a class of  large amplitude solutions toward the global Maxwellian with an explicit almost exponential rate in large time was also obtained by   Desvillettes and Villani \cite{Desvillettes-V} conditionally under  some assumptions on smoothness and polynomial moment bounds of the solutions. The result has been recently improved by Gualdani, Mischler and Mouhot \cite{GMM} to derive a sharp exponential time rate by developing an abstract semigroup theory for linear operators which are non-symmetric in some Banach spaces.

On the other hand, in the perturbation framework, i.e., for the case when the solution is sufficiently close to a global Maxwellian  in some sense, due to the extensive study of the linearized operator (Grad \cite{Gra}, Ellis and Pinsky \cite{EP}, and Baranger and Mouhot \cite{BM}, for instance), the well-posedness theory of the Boltzmann equation is indeed well established in different kinds of settings since the pioneering work by Ukai \cite{Uk}. For instance, the energy method in smooth Sobolev spaces was developed in Guo \cite{Guo-03} and Liu, Yang and Yu \cite{Liu-Yang-Yu}. Another $L^2\cap L^\infty$ approach was found by Guo \cite{Guo2,Guo} even for treating the Boltzmann equation on a general bounded domain. Note that for the hard sphere model in the torus case, a non-symmetric energy method was also developed in \cite{GMM} to obtain the asymptotic stability of solutions to the global Maxwellian with a sharp exponential time rate for initial data $F_0(x,v)$ such that $F_0-\mu$ is small enough in $L^1_vL^\infty_x((1+|v|)^k)$ with some $k>2$; see also a recent work \cite{BG} for the investigation of the Boltzmann equation on the bounded domain in a similar functional setting.    We also refer the interested reader to  \cite{HW} for  the issue of the macroscopic regularity of Boltzmann equation.

We remark that in those works in the perturbation framework mentioned above, initial data are required to have small amplitude around the global Maxwellian.   To our best knowledge, the global existence and uniqueness problem of  solutions to the Boltzmann equation with initial data of large amplitude still remains open. The purpose of this paper is to develop a $L^\infty_xL^1_v\cap L^\infty_{x,v}$ method for the well-posedness theory of the Boltzmann equation when initial data are allowed to have large amplitude. Precisely speaking, we prove the global existence and uniqueness of solutions to the Boltzmann equation in the whole space or torus when 
$$
F_0-\mu\in L^\infty_{x,v} ((1+|v|)^\beta \mu^{-1/2})
$$ 
with some $\beta>\max\{3,3+\gamma\}$ satisfying an additional smallness condition that 
$$
\mathcal{E}(F_0)+\|F_0-\mu\|_{L^1_xL^\infty_v(\mu^{-1/2})}
$$
is small enough. In particular, initial data can have large amplitude oscillations.  Note that the result is valid for the full range of both the soft and  hard potentials, i.e., $-3<\g\leq1$.   Moreover,  in the torus case, we also show that the solutions tend to the global Maxwellian with  exponential convergence rates for the hard potentials and with algebraical rate for the soft potentials.

Now we begin to  formulate the main results of the paper. As in \cite{Guo}, we define a weight function 
\begin{equation}
	w_\b(v):=
	(1+|v|^2)^{\f\beta2},\notag
\end{equation}
and look for solutions in the form
\begin{equation}
\notag
f(t,x,v):=
\f{F(t,x,v)-\mu(v)}{\sqrt{\mu(v)}}.
\end{equation}
The Boltzmann equation \eqref{1.1} is then rewritten as
\begin{align}\label{1.9}
f_t+v\cdot\nabla_xf+Lf=\Gamma(f,f),
\end{align}
where the linearised  term 
is given by
\begin{align}
Lf=\nu(v) f-Kf=-\f1{\sqrt{\mu}}\Big\{Q(\mu,\sqrt{\mu}f)+Q(\sqrt{\mu}f,\mu)\Big\},\notag
\end{align}
with $K:=
K_2-K_1$  defined (cf.~\cite{Desvillettes-V}) as
\begin{align}\notag
(K_1f)(v)&=\int_{\mathbb{R}^3}\int_{\mathbb{S}^2}B(v-u,\t)\sqrt{\mu(v)\mu(u)}f(u)\,d\omega du,
\\
(K_2f)(v)&=\int_{\mathbb{R}^3}\int_{\mathbb{S}^2}B(v-u,\t)\sqrt{\mu(u)\mu(u')}f(v')\,d\omega du,\nonumber\\
&\quad+\int_{\mathbb{R}^3}\int_{\mathbb{S}^2}B(v-u,\t)\sqrt{\mu(u)\mu(v')}f(u')\,d\omega du,
\notag
\\
\nu(v)&=\int_{\mathbb{R}^3}\int_{\mathbb{S}^2}B(v-u,\t)\mu(u)\,d\omega du\sim (1+|v|)^{\g},
\notag
\end{align}
and the nonlinear term is given by
\begin{align}
\Gamma(f,f):=
\f1{\sqrt{\mu}}Q(\sqrt{\mu}f,\sqrt{\mu}f)
&=\f1{\sqrt{\mu}}Q_+(\sqrt{\mu}f,\sqrt{\mu}f)-\f1{\sqrt{\mu}}Q_-(\sqrt{\mu}f,\sqrt{\mu}f)\nonumber\\
&:=
\Gamma_+(f,f)-\Gamma_-(f,f).\notag
\end{align}
Then, from \eqref{1.9}, the mild form of the Boltzmann equation is given by 
\begin{align}\label{1.16}
f(t,x,v)=&e^{-\nu(v)t}f_0(x-vt,v)+\int_0^te^{-\nu(v)(t-s)} (Kf)(s,x-v(t-s),v)ds\nonumber\\
&+\int_0^te^{-\nu(v)(t-s)} \Gamma(f,f)(s,x-v(t-s),v)ds,
\end{align}
for $t\geq 0$, $x\in \Omega$, $v\in \R^3$.

The first result of this paper is stated as follows.

\begin{theorem}[Global Existence]\label{thm1.1}
  Let  $\Omega=\mathbb{T}^3~\mbox{or}~\mathbb{R}^3$. For given $\b>\max\{3,3+\g\}$, $\bar{M}\geq 1$, suppose the initial data $F_0$  satisfies $F_0(x,v)=\mu(v)+\sqrt{\mu(v)}f_0(x,v)\geq 0$ and $\|w_\b f_0\|_{L^\infty}\leq \bar{M}$. Then there is a  small constant $\v_0>0 $  depending on $\g,\b,\bar{M}$ such that if 
  \begin{equation}\label{1.18}
\mathcal{E}(F_0)+\|f_0\|_{L^1_xL^\infty_v}\leq \v_0,
  \end{equation}
the Boltzmann equation \eqref{1.1}, \eqref{1.4}, \eqref{1.5-1} has a global unique  mild solution $F(t,x,v)=\mu(v)+\sqrt{\mu(v)}f(t,x,v)\geq0$  satisfying \eqref{1.13}-\eqref{1.14-1} and 
 \begin{align}
 \|w_\b f(t)\|_{L^\infty}\leq \tilde{C}_1\bar{M}^2,\notag
 \end{align}
where $\tilde{C}_1$ depends only on $\g,\b$. Moreover, if the initial data $f_0$ is continuous in $(x,v)\in\Omega\times\mathbb{R}^3$, then the solution $f(t,x,v)$ is continuous in $[0,\infty)\times\Omega\times\mathbb{R}^3$.
\end{theorem}

\begin{remark}
It should be pointed out that initial data satisfying the smallness condition \eqref{1.18} are  allowed to have  large amplitude oscillations in spatial variable. For instance, one may take
\begin{align}
F_0(x,v)=\rho_0(x)\mu=\f{\r_0(x)}{(2\pi)^{\f32}}e^{-\f{|v|^2}{2}},\quad (x,v)\in\Omega\times\mathbb{R}^3,\notag
\end{align}
with $\r_0(x)\geq0$,  $\r_0\in L^\infty_x$, $\r_0-1\in L^1_x$ and $\r_0\ln\r_0-\r_0+1\in L^1_x$.  Then, it is straightforward to verify that  \eqref{1.18} holds if  $ \|\r_0\ln\r_0-\r_0+1\|_{L^1_x}+\|\r_0-1\|_{L^1_x}$ is small. Even though $ \|\r_0\ln\r_0-\r_0+1\|_{L^1}+\|\r_0-1\|_{L^1_x}$ is required to be small,  initial data are allowed to have large amplitude oscillations. 
\end{remark}

\begin{remark}
From the proof of Theorem \ref{thm1.1} later on, by the same argument, the smallness condition \eqref{1.18} can be relaxed to 
\begin{equation}
\mathcal{E}(F_0)+\sup_{(t,x)\in[t_1,\infty)\times\Omega
}\int_{\mathbb{R}^3}e^{-\nu(v)t}|f_0(x-vt,v)|dv\leq \v_0,\notag
\end{equation}
 where $t_1:=
 (8\tilde{C}_4[1+\|w_\b f_0\|_{L^\infty}])^{-1}$  
 is defined in Proposition \ref{prop7.1} later on.
\end{remark}

\begin{remark}
Under the assumptions of Theorem \ref{thm1.1}, and further let $\b$ suitably large. Let the initial data  $f_0(x,v)\in C^1(\Omega\times\mathbb{R}^3)$ and 
$\|w_{\f{\b}{2}}\nabla_xf_0\|_{L^\infty}+\|w_{\f{\b}{2}}\nabla_vf_0\|_{L^\infty}<+\infty$, then the Boltzmann solution $f(t,x,v)$ obtained in Theorem \ref{thm1.1} satisfies $f(t,x,v)\in C^1(\mathbb{R}_+\times\Omega\times\mathbb{R}^3)$ and 
$$\|w_{\f{\b}{2}}\nabla_xf(t)\|_{L^\infty}+\|w_{\f{\b}{2}}\nabla_vf(t)\|_{L^\infty}\leq  \exp\{C(t)\} $$
where  $C(t)>0$ is a continuous function of $t>0$, and depends only on $\bar{M}$ and $\|w_{\f{\b}{2}}\nabla_xf_0\|_{L^\infty}+\|w_{\f{\b}{2}}\nabla_vf_0\|_{L^\infty}$. It should be pointed  out that the above regularity result can be proved by using similar arguments as in the proof of Proposition \ref{prop7.1} in the appendix and Gronwall inequality.
\end{remark}

It follows immediately from Theorem \ref{thm1.1} that even if initial density $\r_0(x):=\int_{\mathbb{R}^3}F_0(x,v)dv$ contains vacuum, then the macroscopic density function 
$$
\rho(t,x):=\int_{\R^3} F(t,x,v)\,dv
$$ 
must have uniformly  positive lower bound in finite time. Indeed, one has
\begin{corollary}[Positive Lower Bound of Density]\label{cor1.3}
Under the same conditions of Theorem \ref{thm1.1}, there exists a positive time $T_0>0$ such that
 \begin{align}
 \notag
 |\r(t,x)-1|=\Big|\int_{\mathbb{R}^3}[F(t,x,v)-\mu(v)]dv\Big|\leq \f34,
 \end{align}
 for all $t\geq T_0$ and $x\in \Omega$.
\end{corollary}

Moreover, for the global solutions obtained in Theorem \ref{thm1.1} with $\Omega=\mathbb{T}^3$, one can further obtain the explicit rates of convergence of solutions in $L^\infty_{x,v}$. Therefore, it shows that even if initial data $f_0(x,v)$ could be large in $L^\infty_{x,v}$, the solution $f(t,x,v)$ must tend to zero as time goes to infinity. In fact, one has  
\begin{theorem}[Decay Estimate for Hard Potentials]\label{thm1.2}
Let  $\Omega=\mathbb{T}^3$,  $0\leq\g\leq 1$, and $\b>\max\{3, 3+\g\}$.  Assume  $(M_0,J_0,E_0)=(0,0,0)$,  and  $\v_0>0$ sufficiently small, then there exists a positive constant $\sigma_0>0$ such that the  solution $f(t,x,v)$ obtained in Theorem \ref{thm1.1} satisfies
\begin{align}\label{1.29}
\|w_\b f(t)\|_{L^\infty}\leq \tilde{C}_2e^{-\sigma_0t},
\end{align}
for all $t\geq 0$, where $\tilde{C}_2>0$ is a positive constant  depending only on the initial data.
\end{theorem}

\begin{theorem}[Decay Estimate for Soft Potentials]\label{thm1.3}
Let $\Omega=\mathbb{T}^3$, $-3<\g<0$,  and $\b\geq\max\{\f92,4+|\g|\}$. Let $\d$ be any given positive constant such that $\d\in(0,\f13)$.  Assume  $(M_0,J_0,E_0)=(0,0,0)$,  and $\v_0>0$ is sufficiently small, then the solution $f(t,x,v)$ obtained in Theorem \ref{thm1.1} satisfies
\begin{align}\label{1.30}
\| f(t)\|_{L^\infty}\leq \tilde{C}_3(1+t)^{-1-\f2{|\g|}+\d},
\end{align}
for all $t\geq 0$, where $\tilde{C}_3>0$ is a positive constant depending only on the initial data.
\end{theorem}

Now we explain the strategy of the proof of the above main results. As mentioned before,  the  only global existence of large-data solutions to the Boltzmann equation is due to DiPerna and Lions \cite{D-Lion} by the weak compactness argument, but the uniqueness  of these renormalized solutions is completely open due to the lack of $L^\infty$ a priori estimates.  Indeed, it is difficult to establish the global $L^\infty$ bound for the solutions of Boltzmann equations due to the nonlinear term $\Gamma(f,f)(t)$. In those aforementioned references \cite{Guo2,Guo,Kim2014,Strain,Ukai-Yang}, one usually  has to estimate the nonlinear term in the following way
\begin{align}
|w_\b(v)\Gamma(f,f)(t)|\leq C\nu(v)\|w_\b f(t)\|^2_{L^\infty},\notag
\end{align}
so that the smallness  assumption on the $L^\infty$-norm is necessarily required.


To remove the above smallness assumption on the $L^\infty$-norm, we need a new idea to control the nonlinear term $\Gamma(f,f)$. For this, we firstly establish a new estimate for the  nonlinear term (see Lemma \ref{lem4.2} below), i.e., for $\b\geq\f12$,
\begin{align}
&\Big|w_\b(v)\Gamma(f,f)(t,x,v)\Big|
\leq  C\nu(v)\|w_{\b}f(t)\|^{2-a}_{L^\infty}\cdot\Big(\int_{\mathbb{R}^3}|f(t,x,u)|du\Big)^{a},\notag
\end{align}
for some $0<a<1$.  Secondly, under the condition \eqref{1.18}, we observe that $\int_{\mathbb{R}^3}|f(t,x,u)|du$ could be small after some positive time, even if it could be initially large due to the hyperbolicity of the Boltzmann equation.   This observation is the key point of this paper to control the nonlinear term $\Gamma(f,f)$. In such way, through careful analysis one can finally obtain the following uniform estimate
\begin{align}
\sup_{0\leq s\leq t}\|w_\b f(s)\|_{L^\infty}\leq C\bar{M}^2,\notag
\end{align}
under smallness of $\|f(t)\|_{L^\infty_x L^1_v}$ uniformly for all $t\geq t_1$ with some $t_1>0$. 
In the whole proof, we shall use  only the smallness of $\mathcal{E}(F_0)+\|f_0\|_{L^1_xL^\infty_v}$ so that  initial data are allowed to have large amplitude oscillations.

\

The paper is organized as follows. In Section 2,  we introduce the  local existence of solutions to  the Boltzmann equation and list some properties on the kernel of linearized operator, and the detailed  proofs can be found in appendix.  In Section 3, we develop the $L^\infty_xL^1_v\cap L^\infty_{x,v}$ estimate to prove the main Theorem \ref{thm1.1}.  The time-decay estimates of the Boltzmann equation on torus are established in Section 4.  

\

\noindent{\bf Notations.}  Throughout this paper, $C$ denotes a generic positive constant  which may depend on $\g,\b$ and  vary from line to line.  $C_a,C_b,\cdots$ denote the generic positive constants depending on $a,~b,\cdots$, respectively, which also may vary from line to line. $\|\cdot\|_{L^2}$ denotes the standard $L^2(\Omega\times\mathbb{R}^3_v)$-norm, and $\|\cdot\|_{L^\infty}$ denotes the $L^\infty(\Omega\times\mathbb{R}^3_v)$-norm.


\section{Preliminaries}

As mentioned before, Kaniel-Shinbrot \cite{Kaniel-S} investigated the local existence and uniqueness of solutions to the Boltzmann equation  for large initial data around vacuum.  Though, to prove Theorem \ref{thm1.1}, we need to figure out more quantitative properties of the local existence regarding the lifespan of the local $L^\infty$ solution in terms of the $L^\infty$ bound of initial data. Therefore, we would give a representation of the  local existence and uniqueness of solutions to the Boltzmann equation applicable for the global $L^\infty$ estimates in our own setting. The proof of the following result will be given in the appendix.

\begin{proposition}[Local Existence]\label{prop7.1}
	Let $\Omega=\mathbb{T}^3$ or $\mathbb{R}^3$,  $-3<\g\leq 1$, $\b>3$,  $F_0(x,v)=\mu(v)+\sqrt{\mu(v)}f_0(x,v)\geq 0$ and $\|w f_0\|_{L^\infty}<\infty$, then there exists a positive time 
	\begin{align}\label{LT}
	t_1:=
	(8\tilde{C}_4[1+\|w_\b f_0\|_{L^\infty}])^{-1}>0,
	\end{align}
	such that the Boltzmann equation \eqref{1.1}, \eqref{1.4}, \eqref{1.5-1} has a unique solution $F(t,x,v)=\mu(v)+\sqrt{\mu(v)}f(t,x,v)\geq 0$ satisfying 
	\begin{align}
	\notag
	\sup_{0\leq t\leq t_1}\|w_\b f(t)\|_{L^\infty}\leq 2\|w_\b f_0\|_{L^\infty},
	\end{align}
	where the positive constant $\tilde{C}_4\geq1$ depending only on $\g,\b$. Moreover,  the conservations of defect mass, momentum, energy \eqref{1.13}-\eqref{1.13-2} as well as the additional defect entropy inequality \eqref{1.14-1} hold. 
	Finally, if initial data $f_0$ are continuous, then the solution $f(t,x,v)$ is continuous in $[0,t_1]\times\Omega\times\mathbb{R}^3$.
\end{proposition}

For later use, we list the following  result on the operator $K$, whose proof will be given in the appendix.  Interested readers may also refer to  \cite{Glassey,Bellomo} for more details.
\begin{lemma}\label{lem2.1}
	For $-3<\g\leq1$, the following Grad's estimates hold
	\begin{align}
	K_1f(v)=\int_{\mathbb{R}^3}k_1(v,\eta)f(\eta)\,d\eta,\quad
	K_2f(v)=\int_{\mathbb{R}^3}k_2(v,\eta)f(\eta)\,dv,\notag
	\end{align}
	where  $k_1(v,\eta)$ and $k_2(v,\eta)$ satisfy
	\begin{align}
	0\leq k_1(v,\eta)=c_1|v-\eta|^\g e^{-\f{|v|^2}{4}}e^{-\f{|\eta|^2}{4}},\notag
	\end{align}
	and
	\begin{align}\label{2.15}
	0\leq k_2(v,\eta)\leq \f{C_\g}{|v-\eta|^{\f{3-\g}2}}e^{-\f{|v-\eta|^2}{8}}e^{-\f{||v|^2-|\eta|^2|^2}{8|v-\eta|^2}},
	\end{align}
	where $c_1>0$ is a given constant,  and  $C_\g$ is  a constant depending only on $\g$.
	
\end{lemma}

\begin{remark}
Note that the upper bound in \eqref{2.15} is not optimal, but it is enough for the  use of the later proof. Moreover, we will not make any effort on the optimal estimates related to $K$ in order to show Theorem \ref{thm1.1}. 
\end{remark}

From Lemma \ref{lem2.1}, one has that 
\begin{align}
Kf=\int_{\mathbb{R}^3}k(v,\eta)f(\eta)\,d\eta:=
\int_{\mathbb{R}^3}\Big\{k_2(v,\eta)-k_1(v,\eta)\Big\}f(\eta)\,d\eta,\notag
\end{align}
with 
\begin{align}\label{2.17}
|k(v,\eta)|\leq c_1|v-\eta|^\g e^{-\f{|v|^2}{4}}e^{-\f{|\eta|^2}{4}}
+\f{C_\g}{|v-\eta|^{\f{3-\g}2}}e^{-\f{|v-\eta|^2}{8}}e^{-\f{||v|^2-|\eta|^2|^2}{8|v-\eta|^2}}.
\end{align}
By the same calculations as in \cite{Bellomo,Glassey}, it is straightforward to check that for $\a\geq 0$,
\begin{align}\label{2.18}
\int_{\mathbb{R}^3}\Big|k(v,\eta)\cdot\f{w_{\a}(v)}{w_{\a}(\eta)}\Big|\,d\eta\leq C_\g(1+|v|)^{-1}.
\end{align}

In order to deal with difficulties in the case of the soft potentials, as in \cite{Strain-Guo} we introduce  a smooth cutoff function $0\leq\chi_m\leq 1$ with $0<m\leq 1$ such that 
\begin{equation}
\chi_m(s)=1~~\mbox{for}~s\leq m,~~~\chi_m(s)=0~~\mbox{for}~s\geq2m.\notag
\end{equation}
Then we define
\begin{align}
(K^mg)(v)&=\int_{\mathbb{R}^3}\int_{\mathbb{S}^2}B(v-u,\t)\chi_m(|v-u|)\sqrt{\mu(u)\mu(u')}f(v')\,d\omega du\nonumber\\
&\quad+\int_{\mathbb{R}^3}\int_{\mathbb{S}^2}B(v-u,\t)\chi_m(|v-u|)\sqrt{\mu(u)\mu(v')}f(u')\,d\omega du\nonumber\\
&\quad-\int_{\mathbb{R}^3}\int_{\mathbb{S}^2}B(v-u,\t)\chi_m(|v-u|)\sqrt{\mu(v)\mu(u)}f(u)\,d\omega du\nonumber\\
&:=
K_2^mf(v)-K^m_1f(v),\notag
\end{align}
and 
\begin{equation}\label{2.30}
K^c=K-K^m.
\end{equation}

The following result on $K^m$ and $K^c$ can be regarded as a refined version of  \cite[Lemma 1]{Strain-Guo}, and its proof  can be found  in the appendix.

\begin{lemma}\label{lem2.2}
Let  $-3<\gamma\leq1$,  then it holds that 
	\begin{equation}\label{2.31}
	|(K^mg)(v)|\leq Cm^{3+\gamma}e^{-\f{|v|^2}{10}}\|g\|_{L^\infty},
	\end{equation}
	and 
	\begin{equation}
	(K^cg)(v)=\int_{\mathbb{R}^3}l(v,\eta)g(\eta)\,d\eta.\notag
	\end{equation}
	Here the kernel $l(v,\eta)$ satisfies that for $0\leq a\leq 1$,
	\begin{align}\label{2.33}
	|l(v,\eta)|&\leq\f{C_\g m^{a(\g-1)}}{|v-\eta|^{1+\f{(1-a)}{2}(1-\g)}}\f{1}{(1+|v|+|\eta|)^{a(1-\g)}}e^{-\f{|v-\eta|^2}{10}}e^{-\f{||v|^2-|\eta|^2|^2}{16|v-\eta|^2}}\nonumber\\
	&\quad+C|v-\eta|^\g e^{-\f{|v|^2}{4}}e^{-\f{|\eta|^2}{4}},
	\end{align}	
where $C_\g$  is a constant depending only on $\g$. It is worth to  point out that $C_\g$ is uniform in $ a\in[0, 1]$.
\end{lemma}



Since the constant $C_\g$ in \eqref{2.33} does not depend on $a\in[0,1]$,  we have the following estimates on $l(v,\eta)$ from  Lemma \ref{lem2.2} by taking $a=1$ and $a=0$, respectively. 

\begin{lemma}
Let $-3<\g\leq 1$, both the following two bounds on $l(l,v)$ hold:
	\begin{align}\label{2.33-2}
	|l(v,\eta)|&\leq\f{C_\g m^{\g-1}}{|v-\eta|(1+|v|+|\eta|)^{1-\g}}e^{-\f{|v-\eta|^2}{10}}e^{-\f{||v|^2-|\eta|^2|^2}{16|v-\eta|^2}}+C|v-\eta|^\g e^{-\f{|v|^2}{4}}e^{-\f{|\eta|^2}{4}},
	\end{align}	
	and 
	\begin{align}\label{2.33-3}
	|l(v,\eta)|&\leq\f{C_\g}{|v-\eta|^{\f{3-\g}{2}}}e^{-\f{|v-\eta|^2}{10}}e^{-\f{||v|^2-|\eta|^2|^2}{16|v-\eta|^2}}+C|v-\eta|^\g e^{-\f{|v|^2}{4}}e^{-\f{|\eta|^2}{4}}.
	\end{align}		
Moreover,  it holds that 
\begin{equation}\label{2.40}
	\int_{\mathbb{R}^3}\Big|l(v,\eta)\cdot\f{w_{\a}(v)}{w_{\a}(\eta)}\Big|d\eta\leq C_\g m^{\g-1}(1+|v|)^{\g-2}+C_\g e^{-\f{|v|^2}{4}}\leq C(\g)m^{\g-1}\f{\nu(v)}{(1+|v|)^2},
\end{equation}
and 	
\begin{align}\label{2.40-1}
\int_{\mathbb{R}^3}\Big|l(v,\eta)\cdot\f{w_{\a}(v)}{w_{\a}(\eta)}\Big|d\eta&\leq C_\g(1+|v|)^{-1}+C_\g e^{-\f{|v|^2}{4}}\leq C_\g(1+|v|)^{-1},
\end{align}
where $\a\ge0$ is an arbitrary positive constant.
\end{lemma}

\begin{remark}
Indeed, the estimate \eqref{2.33-3} and \eqref{2.40-1} are  the same  as the ones  in \eqref{2.17} and \eqref{2.18}.  On the other hand, the estimates \eqref{2.33-2} and \eqref{2.40} imply that one can get more decay with respect to $v$, but at the cost of  growth with respect to the parameter $\f1{m}$. All these properties will be used later.
\end{remark}


Motivated by Guo \cite{Guo}, we have the following lemma which will be used later.
\begin{lemma}[\cite{Guo}]\label{lem2.5}
	Let $F(t,x,v)$ satisfy \eqref{1.13}, \eqref{1.13-2} and the additional defect entropy inequality \eqref{1.14-1},
	then it holds that 
	\begin{align}\label{2.9}
	&\int_{\Omega}\int_{\mathbb{R}^3}\f{|F(t,x,v)-\mu(v)|^2}{4\mu(v)}I_{\{|F(t,x,v)-\mu(v)|\leq \mu(v)\}}dvdx\nonumber\\
	&\quad+\int_{\Omega}\int_{\mathbb{R}^3}\f{1}{4}|F(t,x,v)-\mu(v)|I_{\{|F(t,x,v)-\mu(v)|\geq \mu(v)\}}dvdx\nonumber\\
	&\leq \int_{\Omega}\int_{\mathbb{R}^3}F_0\ln F_0-\mu\ln\mu dvdx+[\f{3}{2}\ln(2\pi)-1]M_0+\f12E_0=
	\mathcal{E}(F_0).
	\end{align}
\end{lemma}

\noindent{\bf Proof.}  By Taylor expansion, we have 
\begin{equation*}
F(t)\ln F(t)-\mu\ln\mu=(1+\ln\mu)[F(t)-\mu]+\f1{2\tilde{F}}|F(t)-\mu|^2,
\end{equation*}
where $\tilde{F}$ is between $F(t)$ and $\mu$. Noting $1+\ln\mu=-[\f{3}{2}\ln(2\pi)-1]-\f12|v|^2$, we have
\begin{align}\label{2.13}
&\int_{\Omega}\int_{\mathbb{R}^3}\f1{2\tilde{F}}|F(t)-\mu|^2dvdx\nonumber\\
&=\int_{\Omega}\int_{\mathbb{R}^3}[F(t)\ln F(t)-\mu\ln\mu]dvdx
+[\f{3}{2}\ln(2\pi)-1]\int_{\Omega}\int_{\mathbb{R}^3} [F(t)-\mu] dvdx\nonumber\\
&\quad+\f12\int_{\Omega}\int_{\mathbb{R}^3}|v|^2 [F(t)-\mu] dvdx\nonumber\\
&\leq \int_{\Omega}\int_{\mathbb{R}^3}[F_0\ln F_0-\mu\ln\mu]dvdx+[\f{3}{2}\ln(2\pi)-1]M_0+\f12E_0,
\end{align}
where we have used \eqref{1.13}, \eqref{1.13-2}  and \eqref{1.14-1} in the last inequality. We note that $|F-\mu|\geq\mu$ yields that $F\geq2\mu$ or $F=0$, thus we have 
$$\f{|F-\mu|}{\tilde{F}}\geq \f12,$$
which, together with \eqref{2.13}, yields \eqref{2.9}. Therefore, the proof of this lemma is completed. $\hfill\Box$


\section{Global Estimates}
In order to prove the global existence of solutions to the Boltzmann equation, it suffices to get uniform estimates on solutions since one has already obtained in Proposition \ref{prop7.1} the local existence of unique solutions to  the Boltzmann equation with possibly large initial data. In this section, we devote ourselves to establish the global uniform estimate for  the obtained solutions to the Boltzmann equation.

\subsection{Weighted $L^\infty$-Estimate}

Define
\begin{equation}
h(t,x,v):=
w_{\b}(v)f(t,x,v).\notag
\end{equation}
Multiplying \eqref{1.9} by $w_\b(v)$, one gets that 
\begin{align}\label{4.1}
h_t+v\cdot\nabla_xh+\nu(v)h-K_{w_{\b}}h=\Gamma_{w_{\b}}(h,h),
\end{align}
where 
\begin{align}
(K_{w_{\b}}h)= w_{\b}(v)\Big(K\f{h}{w_{\b}}\Big)(v)&=w_{\b}(v)\Big(K^m\f{h}{w_{\b}}\Big)(v)+w_{\b}(v)\Big(K^c\f{h}{w_{\b}}\Big)(v)\nonumber\\
&:=
K_{w_{\b}}^mh+K^c_{w_{\b}}h,\notag
\end{align}
and 
\begin{align}\label{4.3}
\Gamma_{w_{\b}}(h,h)&:=
w_{\b}(v)\Gamma(\f{h}{w_{\b}},\f{h}{w_{\b}})= w_{\b}(v)\Gamma_+(\f{h}{w_{\b}},\f{h}{w_{\b}})- w_{\b}(v)\Gamma_-(\f{h}{w_{\b}},\f{h}{w_{\b}})\nonumber\\
&\equiv w_{\b}(v)\Gamma(f,f)= w_{\b}(v)\Gamma_+(f,f)- w_{\b}(v)\Gamma_-(f,f).
\end{align}
Then the mild solution of \eqref{4.1} can be written as
\begin{align}\label{4.4}
h(t,x,v)=&e^{-\nu(v)t}h_0(x-vt,v)+\int_0^te^{-\nu(v)(t-s)} (K_{w_{\b}}^mh)(s,x-v(t-s),v)\,ds\nonumber\\
&+\int_0^te^{-\nu(v)(t-s)} (K_{w_{\b}}^c h)(s,x-v(t-s),v)\,ds\nonumber\\
&+\int_0^te^{-\nu(v)(t-s)} (\Gamma_{w_{\b}}(h,h))(s,x-v(t-s),v)\,ds.
\end{align}

Firstly, we give 
estimates
on the nonlinear term $\Gamma(f,f)$.
\begin{lemma}\label{lem4.2}
Let $-3<\g\leq 1$.  For $\a\geq0$, it holds that 
\begin{align}\label{4.8-1}
\begin{cases}
\dis |w_{\a}(v)\Gamma_-(f,f)(s,y,v)|
\leq C_\g\nu(v)\|w_{\a}f(s)\|_{L^\infty}\cdot\|f(s)\|^{\f{4p+1}{5p}}_{L^\infty}\cdot\Big(\int_{\mathbb{R}^3}|f(s,y,u)|\,du\Big)^{\f{p-1}{5p}},\\[3mm]
|w_{\a}(v)\Gamma_+(f,f)(s,y,v)|
\dis \leq C_\g\nu(v)\|w_{\a}f(s)\|_{L^\infty}\cdot\|w_{\f12}f(s)\|^{\f{4p+1}{5p}}_{L^\infty}\cdot\Big(\int_{\mathbb{R}^3}|f(s,y,u)|\,du\Big)^{\f{p-1}{5p}},
\end{cases}
\end{align}
where  $p>1$  is defined in \eqref{4.10}.
\end{lemma}

\noindent{\bf Proof.}  
It is noted that 
\begin{align}\label{4.9}
\Big|w_{\a}(v)\Gamma_-(f,f)(s,y,v)\Big|
&\leq C\|w_{\a}f(s)\|_{L^\infty} \int_{\mathbb{R}^3}|v-u|^{\g}\sqrt{\mu(u)}|f(s,y,u)|\,du.
\end{align}
To estimate the integration term on the RHS of \eqref{4.9}, we choose
\begin{equation}\label{4.10}
p=1+\f{3+\g}{4(9-\g)}~~\mbox{for}~-3<\g\leq 1,
\end{equation}
which yields  that
\begin{align}\label{4.11}
1<p\leq {\frac{9}{8}},\quad p_\ast:=
\f{p}{p-1}\geq {9},\quad p\f{\g-3}{2}>-3,\quad p\g>-3.
\end{align}
Then it follows from  \eqref{4.9}, \eqref{4.11} and H\"older inequality that 
\begin{align}\label{4.13}
&\Big|w_{\a}(v)\Gamma_-(f,f)(s,y,v)\Big|\nonumber\\
&\leq C\|w_{\a}f(s)\|_{L^\infty} \Big(\int_{\mathbb{R}^3}|v-u|^{p\g}\sqrt{\mu(u)}du\Big)^{\f1p}
{\Big(\int_{\mathbb{R}^3}\sqrt{\mu(u)}|f(s,y,u)|^{\f{p}{p-1}}du\Big)^{1-\f1p}}\nonumber\\
&\leq C_\g\nu(v)\|w_{\a}f(s)\|_{L^\infty} \Big(\int_{\mathbb{R}^3}|f(s,y,u)|^{\f{5p}{p-1}}du\Big)^{\f{p-1}{5p}}\nonumber\\
&\leq C_\g\nu(v)\|w_{\a}f(s)\|_{L^\infty} {\|f(s)\|^{\f{4p+1}{5p}}_{L^\infty}\Big(\int_{\mathbb{R}^3}|f(s,y,u)|du\Big)^{\f{p-1}{5p}}}.
\end{align}

Next, we consider the gain term which needs much more care.  We note that 
$$|v|^2\leq |u'|^2+|v'|^2,$$ 
which yields
\begin{equation}
\mbox{either}~~\f12|v|^2\leq |u'|^2~~\mbox{or}~~ \f12|v|^2\leq |v'|^2.\notag
\end{equation} 
Hence one obtains that 
\begin{align}\label{4.15}
&\Big|w_{\a}(v)\Gamma_+(f,f)(s,y,v)\Big|=\Big|\f{w_{\a}(v)}{\sqrt{\mu(v)}}Q_+(\sqrt{\mu}f,\sqrt{\mu}f)(s,y,v)\Big|\nonumber\\
&\leq w_\a (v) \int_{\mathbb{R}^3}\int_{\mathbb{S}^2}B(v-u,\t)\sqrt{\mu(u)}\Big|f(s,y,u') f(s,y,v')
\Big| I_{\{\f12|v|^2\leq |u'|^2\}}dud\omega\nonumber\\
&\quad +w_\a(v)\int_{\mathbb{R}^3}\int_{\mathbb{S}^2}B(v-u,\t)\sqrt{\mu(u)}\Big|f(s,y,u') f(s,y,v')
\Big|I_{\{\f12|v|^2\leq |v'|^2\}}dud\omega\nonumber\\
&\leq C \int_{\mathbb{R}^3}\int_{\mathbb{S}^2}B(v-u,\t)\sqrt{\mu(u)}\Big|w_{\a}(u')f(s,y,u')f(s,y,v')
\Big| dud\omega \nonumber\\
&\quad +C \int_{\mathbb{R}^3}\int_{\mathbb{S}^2}B(v-u,\t)\sqrt{\mu(u)}\Big|f(s,y,u')w_{\a}(v')f(s,y,v')
\Big|dud\omega\nonumber\\
&=I_{1}+I_{2}.
\end{align}
To estimate $I_{1}$, as in \cite{Glassey,Bellomo},  we use the change of variables
\begin{align}\label{3.9}
z=u-v,~~z_{\shortparallel}=[z\cdot\omega]\omega,~~z_{\perp}=z-z_{\shortparallel},~~ \eta=v+z_{\shortparallel}.
\end{align}
Then it holds that 
\begin{align}\label{3.10}
u'=v+z_{\perp},~~~v'=v+z_{\shortparallel},
\end{align}
and 
\begin{equation}\label{3.11}
B(v-u,\t)\leq C|z_{\shortparallel}|(|z_{\shortparallel}|+|z_{\perp}|)^{\g-1}.
\end{equation}
Hence it follows from \eqref{4.10}, \eqref{4.11},  \eqref{3.9} and \eqref{3.10} that
\begin{align}
I_{1}&=C\int_{\mathbb{R}^3}\int_{\mathbb{S}^2}B(v-u,\t)\sqrt{\mu(u)}\Big|w_\a(u')f(s,y,u')f(s,y,v')
\Big| dud\omega \nonumber\\
&\leq C\|w_{\a}f(s)\|_{L^\infty}\Big(\int_{\mathbb{R}^3}|v-u|^{p\g}\sqrt{\mu(u)}du\Big)^{\f1p}
\Big(\int_{\mathbb{R}^3}\int_{\mathbb{S}^2}\sqrt{\mu(u)}|f(s,y,v')
|^{\f{p}{p-1}}dud\omega\Big)^{1-\f1p} \nonumber\\
&\leq C_\g\nu(v)\|w_{\a}f(s)\|_{L^\infty}\Big(\int_{\mathbb{R}^3}\int_{\mathbb{S}^2}e^{-\f{|v+z_{\shortparallel}+z_{\perp}|^2}{4}}|f(s,y,v+z_{\shortparallel})|^{\f{p}{p-1}}dzd\omega\Big)^{1-\f1p} \nonumber\\
&\leq C_\g\nu(v)\|w_{\a}f(s)\|_{L^\infty}{\Big(\int_{\mathbb{R}^3}\int_{z_{\perp}}\f1{|\eta-v|^2}e^{-\f{|\eta+z_{\perp}|^2}{4}}|f(s,y,\eta)|^{\f{p}{p-1}}dz_{\perp}d\eta\Big)^{1-\f1p}} \nonumber\\
&\leq C_\g\nu(v)\|w_{\a}f(s)\|_{L^\infty}\Big(\int_{\mathbb{R}^3}\f{(1+|\eta|)^{-4}}{|\eta-v|^{\f52}}d\eta\Big)^{\f45(1-\f1p)}\Big(\int_{\mathbb{R}^3}(1+|\eta|)^{16}|f(s,y,\eta)|^{\f{5p}{p-1}}d\eta\Big)^{\f{p-1}{5p}}\nonumber\\
&\leq C_\g\nu(v)\|w_{\a}f(s)\|_{L^\infty}\|(1+|\eta|)^{16\f{p-1}{4p+1}}f(s)\|^{\f{4p+1}{5p}}_{L^\infty}\Big(\int_{\mathbb{R}^3}|f(s,y,\eta)|d\eta\Big)^{\f{p-1}{5p}}\nonumber\\
&\leq C_\g\nu(v)\|w_{\a}f(s)\|_{L^\infty}{\|w_{\f12}f(s)\|^{\f{4p+1}{5p}}_{L^\infty}\Big(\int_{\mathbb{R}^3}|f(s,y,\eta)|d\eta\Big)^{\f{p-1}{5p}}},\notag
\end{align}
where in the last inequality, we have used the fact that $\f{16(p-1)}{4p+1}\leq \f12$. 
On the other hand, it is noted that by a rotation, one obtains the interchange of $v'$ and $u'$, and then $I_{2}$ can be   changed to  a form similar to $I_{1}$.  Hence, for $I_{2}$, one can also obtain the same estimate as above. Thus one can get that 
\begin{align}\label{4.17}
&\Big|w_\a(v)\Gamma_+(f,f)(s,y,v)\Big|
\leq C_\g\nu(v)\|w_{\a}f(s)\|_{L^\infty}\|w_{\f12}f(s)\|^{\f{4p+1}{5p}}_{L^\infty}\Big(\int_{\mathbb{R}^3}|f(s,y,\eta)|d\eta\Big)^{\f{p-1}{5p}}.
\end{align}
Then \eqref{4.8-1} follows from \eqref{4.13} and \eqref{4.17}.  The proof of Lemma \ref{lem4.2} is completed.   $\hfill\Box$

\

\begin{lemma}\label{lem4.1}
Let $\b> 3$ and  $-3<\g\leq 1$, then it holds that 
\begin{align}\label{4.5}
\sup_{0\leq s\leq t}\|h(s)\|_{L^\infty}
&\leq C_1\Big\{\|h_0\|_{L^\infty}+\|h_0\|^2_{L^\infty}+\sqrt{\mathcal{E}(F_0)}+\mathcal{E}(F_0)\Big\}\nonumber\\
&\quad+C_1\sup_{t_1\leq s\leq t,~y\in\Omega}\left\{ \|h(s)\|^{\f{9p+1}{5p}}_{L^\infty}\Big(\int_{\mathbb{R}^3}|f(s,y,\eta)|d\eta\Big)^{\f{p-1}{5p}}\right\},
\end{align}
where the positive constant $C_1\geq1$ depends only on $\g,\b$, the lifespan  $t_1>0$ is defined in \eqref{LT}, and $p>1$ is defined in \eqref{4.10}.
\end{lemma}

\noindent{\bf Proof.} It follows from \eqref{4.4} that
\begin{align}
|h(t,x,v)|&\leq e^{-\nu(v)t}\|h_0\|_{L^\infty}+\int_0^te^{-\nu(v)(t-s)} \Big|(K_{w_{\b}}^mh)(s,x-v(t-s),v)\Big|ds\nonumber\\
&\quad+\int_0^te^{-\nu(v)(t-s)} \Big|(K_{w_{\b}}^c h)(s,x-v(t-s),v)\Big|ds\nonumber\\
&\quad+\int_0^te^{-\nu(v)(t-s)} \Big|(\Gamma_{w_{\b}}(h,h))(s,x-v(t-s),v)\Big|ds\nonumber\\
&=e^{-\nu(v)t}\|h_0\|_{L^\infty}+J_1+J_2+J_3.\notag
\end{align}
Using \eqref{2.31}, one gets  that 
\begin{align}
J_1&=
\int_0^te^{-\nu(v)(t-s)} \Big|w_{\b}(v)\Big(K^mf\Big)(s,x-v(t-s),v)\Big|ds\nonumber\\
&\leq Cm^{3+\g}\int_0^te^{-\nu(v)(t-s)} w_{\b}(v)e^{-\f{|v|^2}{10}}\|f(s)\|_{L^\infty}ds\nonumber\\
&\leq Cm^{3+\g}e^{-\f{|v|^2}{20}}\sup_{0\leq s\leq t}\|f(s)\|_{L^\infty}.\label{4.7-1}
\end{align}
For $J_3$,   it follows from \eqref{4.3} and \eqref{4.8-1} that for $\b\geq 1/2$,
\begin{align}\label{4.19}
J_3&=\int_0^te^{-\nu(v)(t-s)} |\Gamma_{w_\b}(h,h)(s,x-v(t-s),v)|ds\nonumber\\
&\leq C\sup_{0\leq s\leq t,~y\in\Omega}\left\{ \|h(s)\|^{\f{9p+1}{5p}}_{L^\infty}\Big(\int_{\mathbb{R}^3}|f(s,y,\eta)|d\eta\Big)^{\f{p-1}{5p}}\right\}.
\end{align}

Let $l_{w_\b}(v,v')$ be the corresponding kernel associated with $K^c_{w_\b}$, then we have that
\begin{align}
l_{w_\b}(v,v')=l(v,v')\f{w_\b(v)}{w_\b(v')},\notag
\end{align}
which together with \eqref{2.40} and \eqref{2.40-1}, yield  that 
\begin{equation}\label{4.20-1}
\int_{\mathbb{R}^3}\Big|l_{w_\b}(v,v')\Big| dv'\leq C_\g m^{\g-1}\f{\nu(v)}{(1+|v|)^2}~\mbox{and}~
\int_{\mathbb{R}^3}\Big|l_{w_\b}(v,v')\Big| dv'\leq C_\g (1+|v|)^{-1}.
\end{equation}
For $J_2$, denoting $\tilde{x}:=
x-v(t-s)$, we note that 
\begin{align}
J_2\leq \int_0^te^{-\nu(v)(t-s)} \int_{\mathbb{R}^3}|l_{w_\b}(v,v') h(s,\tilde{x},v')|dv'ds.\notag
\end{align}
To bound the above term, similar as in \cite{Vidav,Guo}, we use \eqref{4.4} again to get that 
\begin{align}\label{4.23}
J_2&\leq \|h_0\|_{L^\infty}\int_0^te^{-\nu(v)(t-s)} \int_{\mathbb{R}^3}|l_{w_\b}(v,v')|dv' ds\nonumber\\
&\quad+\int_0^te^{-\nu(v)(t-s)} \int_{\mathbb{R}^3}|l_{w_\b}(v,v')| \int_0^se^{-\nu(v')(s-\tau)} |(K_{w_\b}^mh)(\tau,\tilde{x}-v'(s-\tau),v')|d\tau dv'ds\nonumber\\
&\quad+\int_0^te^{-\nu(v)(t-s)} \int_{\mathbb{R}^3}|l_{w_\b}(v,v')| \int_0^se^{-\nu(v')(s-\tau)} |(\Gamma_{w_\b}(h,h))(\tau,\tilde{x}-v'(s-\tau),v')|d\tau dv'ds\nonumber\\
&\quad+\int_0^te^{-\nu(v)(t-s)} \int_{\mathbb{R}^3}\int_{\mathbb{R}^3}|l_{w_\b}(v,v')l_{w_\b}(v',v'')|\nonumber\\
&\qquad\qquad\qquad\qquad\qquad\times \int_0^se^{-\nu(v')(s-\tau)} |h(\tau,\tilde{x}-v'(s-\tau),v'')|d\tau dv'dv''ds \nonumber\\
&:=
J_{21}+J_{22}+J_{23}+J_{24}.
\end{align}
It follows from  \eqref{4.19}  and \eqref{4.20-1} that for $\b\geq\f12$,  
\begin{align}
 J_{21}+J_{23}
&\leq C\left\{\|h_0\|_{L^\infty}+\sup_{0\leq s\leq t,~y\in\Omega}\Big( \|h(s)\|^{\f{9p+1}{5p}}_{L^\infty}\Big(\int_{\mathbb{R}^3}|f(s,y,\eta)|d\eta\Big)^{\f{p-1}{5p}}\Big)\right\} \nonumber\\
&\quad\times \int_0^te^{-\nu(v)(t-s)} \int_{\mathbb{R}^3}|l_{w_\b}(v,v')| dv'ds\nonumber\\
&\leq Cm^{\g-1}\left\{\|h_0\|_{L^\infty}+\sup_{0\leq s\leq t,~y\in\Omega}\Big( \|h(s)\|^{\f{9p+1}{5p}}_{L^\infty}\Big(\int_{\mathbb{R}^3}|f(s,y,\eta)|d\eta\Big)^{\f{p-1}{5p}}\Big)\right\}.\notag
\end{align}
For $J_{22}$,  using \eqref{4.7-1},  one has that
\begin{align}
J_{22}&\leq Cm^{3+\g}\sup_{0\leq s\leq t}\|f(s)\|_{L^\infty}\int_0^te^{-\nu(v)(t-s)} \int_{\mathbb{R}^3}|l_{w_\b}(v,v')|e^{-\f{1}{20}|v'|^2}dv'ds\nonumber\\
&\leq Cm^{3+\g}\sup_{0\leq s\leq t}\|f(s)\|_{L^\infty},\notag
\end{align}
where we have used the fact that
\begin{align}
\int_{\mathbb{R}^3}|l_{w_\b}(v,v')|e^{-\f{1}{20}|v'|^2}dv'\leq Ce^{-\f{|v|^2}{100}},\notag
\end{align}
which follows from \eqref{2.33-3} and  similar arguments as in \cite{Glassey}.

\

We now concentrate on the last term $J_{24}$ on the RHS of \eqref{4.23}. As in  \cite{Guo}, we divide it into the following several cases.

\

\noindent{\it Case 1}. For $|v|\geq N$, it follows from \eqref{4.20-1} that
\begin{align}
\int_{\mathbb{R}^3}\Big|l_{w_\b}(v,v')\Big| dv'\leq C_m\f{\nu(v)}{N^2},\notag
\end{align}
which yields immediately that  
\begin{align}
J_{24}&\leq C_m\sup_{0\leq s\leq t}\|h(s)\|_{L^\infty}\int_0^te^{-\nu(v)(t-s)} \int_{\mathbb{R}^3}|l_{w_\b}(v,v')|\int_0^se^{-\nu(v')(s-\tau)} \f{\nu(v')}{1+|v'|^2}d\tau dv'ds\nonumber\\
&\leq \f{C_m}{N^2}\sup_{0\leq s\leq t}\|h(s)\|_{L^\infty}.\notag
\end{align}

\noindent{\it Case 2}. For either $|v|\leq N,~|v'|\geq2N$ or $|v'|\leq2N,~|v''|\geq 3N$. It is noted that we have either $|v-v'|\geq N$ or $|v'-v''|\geq N$, and either one of the following is valid
\begin{align} 
|l_{w_\b}(v,v')|&\leq Ce^{-\f{N^2}{20}}\Big|l_{w_\b}(v,v')e^{\f{|v-v'|^2}{20}}\Big|,\label{4.29}\\
|l_{w_\b}(v',v'')|&\leq Ce^{-\f{N^2}{20}}\Big|l_{w_\b}(v',v'')e^{\f{|v'-v''|^2}{20}}\Big|.
\notag
\end{align}
From \eqref{2.33-2}, a direct calculation shows that
\begin{align}
\small{\int_{\mathbb{R}^3}\Big|l_{w_\b}(v,v')e^{\f{|v-v'|^2}{20}}\Big|dv'\leq C_m\f{\nu(v)}{(1+|v|)^2}~\mbox{and}~\int_{\mathbb{R}^3}\Big|l_{w_\b}(v',v'')e^{\f{|v'-v''|^2}{20}}\Big|dv''\leq C_m\f{\nu(v')}{(1+|v'|)^2}.}\label{4.31}
\end{align}
Then it follows from \eqref{4.29}-\eqref{4.31} that
\begin{align}
&\int_0^te^{-\nu(v)(t-s)} \Big\{\int_{|v|\leq N,|v'|\geq2N}+\int_{|v'|\leq 2N,|v''|\geq3N}\Big\}|l_{w_\b}(v,v')l_{w_\b}(v',v'')|\nonumber\\
&\qquad\times \int_0^se^{-\nu(v')(s-\tau)} |h(\tau,\tilde{x}-v'(s-\tau),v'')|d\tau dv'dv''ds\nonumber\\
&\leq C_me^{-\f{N^2}{20}}\sup_{0\leq s\leq t}\|h(s)\|_{L^\infty}.\notag
\end{align}

\noindent{\it Case 3}. $|v|\leq N,~|v'|\leq 2N,~|v''|\leq 3N$. This is the last remaining case. It is noted that 
\begin{align}\label{4.33}
&\int_0^te^{-\nu(v)(t-s)} \int_{|v'|\leq 2N,|v''|\leq3N}|l_{w_\b}(v,v')l_{w_\b}(v',v'')|\nonumber\\
&\qquad\qquad\times \int_0^se^{-\nu(v')(s-\tau)} |h(\tau,\tilde{x}-v'(s-\tau),v'')|d\tau dv'dv''ds\nonumber\\
&\leq \int_0^te^{-\nu(v)(t-s)} \int_{|v'|\leq 2N,|v''|\leq3N}|l_{w_\b}(v,v')l_{w_\b}(v',v'')|\nonumber\\
&\qquad\qquad\times \int_{s-\l}^se^{-\nu(v')(s-\tau)} |h(\tau,\tilde{x}-v'(s-\tau),v'')|d\tau dv'dv''ds\nonumber\\
&\quad+\int_0^te^{-\nu(v)(t-s)} \int_{|v'|\leq 2N,|v''|\leq3N}|l_{w_\b}(v,v')l_{w_\b}(v',v'')|\nonumber\\
&\qquad\qquad\times \int_0^{s-\l}e^{-\nu(v')(s-\tau)} |h(\tau,\tilde{x}-v'(s-\tau),v'')|d\tau dv'dv''ds.
\end{align}
Using \eqref{4.20-1}, we can bound the first term on the RHS of \eqref{4.33} by
\begin{align}\label{4.34}
C_m\l\sup_{0\leq s\leq t}\|h(s)\|_{L^\infty}\int_0^te^{-\nu(v)(t-s)} \nu(v)ds\leq C_m\l\sup_{0\leq s\leq t}\|h(s)\|_{L^\infty}.
\end{align}

Now we estimate the second term on the RHS of \eqref{4.33}.  Since $l_{w_\b}(v,v')$ has a possible singularity of $\f{1}{|v-v'|}$, we  choose a smooth compact support function $l_N(v,v')$ such that
\begin{align}\label{4.35}
\sup_{|p|\leq 3N}\int_{|v'|\leq 3N}\Big|l_{w_\b}(p,v')-l_N(p,v')\Big|dv'\leq \f{C_m}{N^7}.
\end{align}
Noting
\begin{align}\label{4.36}
l_{w_\b}(v,v')l_{w_\b}(v',v'')&=\Big(l_{w_\b}(v,v')-l_N(v,v')\Big)l_{w_\b}(v',v'')\nonumber\\
&\quad +\Big(l_{w_\b}(v',v'')-l_N(v',v'')\Big)l_N(v,v')+l_N(v,v')l_N(v',v''),
\end{align}
and then using \eqref{4.35} and \eqref{4.36}, we can bound the second term on the RHS of \eqref{4.33} by 
\begin{align}\label{4.37}
&\f{C_m}{N^7}\sup_{0\leq s\leq t}\|h(s)\|_{L^\infty}\int_0^te^{-\nu(v)(t-s)}\int_0^{s-\l}e^{-\nu(v')(s-\tau)} d\tau ds \nonumber\\
&\qquad\qquad\times \left\{\sup_{|v'|\leq2N}\int_{|v''|\leq3N}|l_{w_\b}(v',v'')|dv''+\sup_{|v|\leq N}\int_{|v'|\leq2N}|l_N(v,v')|dv'\right\}\nonumber\\
&+\int_0^te^{-\nu(v)(t-s)} \int_{|v'|\leq 2N,|v''|\leq3N}|l_N(v,v')l_N(v',v'')|\nonumber\\
&\qquad\qquad\times \int_0^{s-\l}e^{-\nu(v')(s-\tau)} |h(\tau,\tilde{x}-v'(s-\tau),v'')|d\tau dv'dv''ds\nonumber\\
&\leq \f{C_m}{N}\sup_{0\leq s\leq t}\|h(s)\|_{L^\infty}\notag\\
&\quad+C_{N,m}\int_0^t\int_0^{s-\l}e^{-c_N(t-s)}e^{-c_N(s-\tau)} 
\int_{|v'|\leq 2N,|v''|\leq3N} |h(\tau,\tilde{x}-v'(s-\tau),v'')|dv'dv''d\tau ds,
\end{align}
where we have used the facts that   $l_N(v,v')l_N(v',v'')$ is bounded and
\begin{equation}
\nu(v)\geq c_N~~\mbox{for}~|v|\leq N,~~\mbox{and}~~\nu(v')\geq c_N~~\mbox{for}~|v'|\leq 2N.\notag
\end{equation}
It follows from \eqref{2.9} and H\"older inequality that 
\begin{align}\label{4.39}
&\int_{|v'|\leq 2N,|v''|\leq3N} |h(\tau,\tilde{x}-v'(s-\tau),v'')|dv'dv''\nonumber\\
&=C_{N,m}\int_{|v'|\leq 2N,|v''|\leq3N} \f{|F(\tau,\tilde{x}-v'(s-\tau),v'')-\mu(v'')|}{\sqrt{\mu(v'')}} I_{\{|F(\tau,\tilde{x}-v'(s-\tau),v'')-\mu(v'')|\leq\mu(v'')\}}dv'dv''\nonumber\\
&\quad+C_{N,m}\int_{|v'|\leq 2N,|v''|\leq3N} |F(\tau,\tilde{x}-v'(s-\tau),v'')-\mu(v'')| I_{\{|F(\tau,\tilde{x}-v'(s-\tau),v'')-\mu(v'')|\geq\mu(v'')\}}dv'dv''\nonumber\\
&\leq C_{N,m}\f{1+(s-\tau)^{\f32}}{(s-\tau)^{\f32}}\left\{\int_{\Omega}\int_{|v''|\leq3N} \f{|F(\tau,y,v'')-\mu(v'')|^2}{\mu(v'')} I_{\{|F(\tau,y,v'')-\mu(v'')|\leq\mu(v'')\}}dv''dy\right\}^{\f12}\nonumber\\
&\quad+C_{N,m}\f{1+(s-\tau)^{3}}{(s-\tau)^{3}}\int_{\Omega}\int_{|v''|\leq3N} |F(\tau,y,v'')-\mu(v'')| I_{\{|F(\tau,y,v'')-\mu(v'')|\geq\mu(v'')\}}dv''dy\nonumber\\
&\leq C_{N,m}\f{1+(s-\tau)^{\f32}}{(s-\tau)^{\f32}}\sqrt{\mathcal{E}(F_0)}+C_{N,m}\f{1+(s-\tau)^{3}}{(s-\tau)^{3}}\mathcal{E}(F_0),
\end{align}
where we have made a change of variable $y=\tilde{x}-v'(s-\tau)$. From \eqref{4.39}, we can bound the second term on the RHS of \eqref{4.37} as follows:
\begin{align}\label{4.40}
&C_{N,m}\int_0^t\int_0^{s-\l}e^{-c_N(t-s)}e^{-c_N(s-\tau)}\int_{|v'|\leq 2N,|v''|\leq3N} |h(\tau,\tilde{x}-v'(s-\tau),v'')|dv'dv''d\tau ds\nonumber\\
&\leq C_{N,m}\l^{-\f32}\sqrt{\mathcal{E}(F_0)}+C_{N,m}\l^{-3}\mathcal{E}(F_0).
\end{align}
Combining \eqref{4.40}, \eqref{4.37}, \eqref{4.34} and \eqref{4.33}, one gets that 
\begin{align}
&\int_0^te^{-\nu(v)(t-s)} \int_{|v'|\leq 2N,|v''|\leq3N}|l_{w_\b}(v,v')l_{w_\b}(v',v'')|\nonumber\\
&\qquad\times \int_0^se^{-\nu(v')(s-\tau)} |h(\tau,\tilde{x}-v'(s-\tau),v'')|d\tau dv'dv''ds\nonumber\\
&\leq C_m(\l+\f{1}{N})\sup_{0\leq s\leq t}\|h(s)\|_{L^\infty}+C_{N,m}\l^{-\f32}\sqrt{\mathcal{E}(F_0)}+C_{N,m}\l^{-3}\mathcal{E}(F_0).\notag
\end{align}

Therefore, collecting  all the above estimates,  we have established that for any $\l>0$ and large $N\geq1$, 
\begin{align}
\sup_{0\leq s\leq t}\|h(s)\|_{L^\infty}&\leq C_m\|h_0\|_{L^\infty}+C\Big(m^{3+\g}+C_m\l+\f{C_m}{N}\Big)\cdot\sup_{0\leq s\leq t}\|h(s)\|_{L^\infty}\nonumber\\
&\quad+C_m\sup_{0\leq s\leq t,~y\in\Omega}\left\{ \|h(s)\|^{\f{9p+1}{5p}}_{L^\infty}\Big(\int_{\mathbb{R}^3}|f(s,y,\eta)|d\eta\Big)^{\f{p-1}{5p}}\right\}\nonumber\\
&\quad+C_{N,m}\l^{-\f32}\sqrt{\mathcal{E}(F_0)}+C_{N,m}\l^{-3}\mathcal{E}(F_0).\notag
\end{align}
Noting $3+\g>0$, first choosing $m$ small, then $\l$ small, and finally letting $N$ be sufficiently large so that $C\Big(m^{3+\g}+C_m\l+\f{C_m}{N}\Big)\leq \f12$, one obtains that for $\b\geq\f12$,
 \begin{align}\label{4.43}
\sup_{0\leq s\leq t}\|h(s)\|_{L^\infty}&\leq C\Big\{\|h_0\|_{L^\infty}+\sqrt{\mathcal{E}(F_0)}+\mathcal{E}(F_0)\Big\}\nonumber\\
&\quad+C\sup_{0\leq s\leq t,~y\in\Omega}\left\{ \|h(s)\|^{\f{9p+1}{5p}}_{L^\infty}\Big(\int_{\mathbb{R}^3}|f(s,y,\eta)|d\eta\Big)^{\f{p-1}{5p}}\right\}.
 \end{align}
Using Proposition \ref{prop7.1}, one has that  for $\b>3$, 
\begin{equation}\label{4.43-1}
\sup_{0\leq s\leq t_1,~y\in\Omega}\left\{ \|h(s)\|^{\f{9p+1}{5p}}_{L^\infty}\Big(\int_{\mathbb{R}^3}|f(s,y,\eta)|d\eta\Big)^{\f{p-1}{5p}}\right\}\leq C\sup_{0\leq s\leq t_1}\|h(s)\|^2_{L^\infty}\leq C\|h_0\|^2_{L^\infty}.
\end{equation}
Substituting \eqref{4.43-1} into \eqref{4.43}, one gets that for $\b>3$, 
 \begin{align}
\sup_{0\leq s\leq t}\|h(s)\|_{L^\infty}&\leq C_1\Big\{\|h_0\|_{L^\infty}+\|h_0\|^2_{L^\infty}+\sqrt{\mathcal{E}(F_0)}+\mathcal{E}(F_0)\Big\}\nonumber\\
 &\quad+C_1\sup_{t_1\leq s\leq t,~y\in\Omega}\left\{ \|h(s)\|^{\f{9p+1}{5p}}_{L^\infty}\Big(\int_{\mathbb{R}^3}|f(s,y,\eta)|d\eta\Big)^{\f{p-1}{5p}}\right\},\notag
 \end{align}
which yields immediately \eqref{4.5}, where the positive constant $C_1$ depends only on $\g,\b$.  Thus the proof of Lemma \ref{lem4.1} is completed. $\hfill\Box$

\subsection{$L^\infty_xL^1_v$ Estimate}

In this subsection, we will concentrate on the estimate of 
\begin{equation}
\int_{\mathbb{R}^3}|f(t,x,v)|dv.\notag
\end{equation}
If $\mathcal{E}(F_0)+\|f_0\|_{L^1_xL^\infty}$ is small, due to the hyperbolicity of the Boltzmann equation one should be able to show that  $\int_{\mathbb{R}^3}|f(t,x,v)|dv$ is small for $t\geq t_1$, even though it could be initially large, i.e.,  $\int_{\mathbb{R}^3}|f_0(x,v)|dv$ is large.  Indeed, we have the following lemma which plays a key role in this paper.
\begin{lemma}\label{lem5.1}
Let  $-3<\g\leq 1$ and $\b>\max\{3,3+\g\}$,  then it holds that
\begin{align}\label{5.2}
\int_{\mathbb{R}^3}|f(t,x,v)|dv&\leq \int_{\mathbb{R}^3}e^{-\nu(v)t}|f_0(x-vt,v)|dv+C_N\l^{-\f32}\sqrt{\mathcal{E}(F_0)}+C_N\l^{-3}\mathcal{E}(F_0)\nonumber\\
&\quad+C\Big(m^{3+\g}+C_m[\l+\f{1}{N}+\f{1}{N^{\b-3}}]\Big)\cdot\sup_{0\leq s\leq t}\Big\{\|h(s)\|_{L^\infty}+\|h(s)\|^2_{L^\infty}\Big\}\nonumber\\
&\quad+C_N\l^{-3}\Big(\sqrt{\mathcal{E}(F_0)}+\mathcal{E}(F_0)\Big)^{1-\f1p}\cdot\sup_{0\leq s\leq t}\|h(s)\|^{1+\f1p}_{L^\infty},
\end{align}
where $\l>0,m>0$ and  $N\geq1$ are to be chosen later. Recall that $p>1$ is defined in \eqref{4.10}.
\end{lemma}

\noindent{\bf Proof.} It follows from \eqref{1.16} and \eqref{2.30} that 
\begin{align}\label{5.3}
&\int_{\mathbb{R}^3}|f(t,x,v)|dv\nonumber\\
&\leq \int_{\mathbb{R}^3}e^{-\nu(v)t}|f_0(x-vt,v)|dv+\int_0^t\int_{\mathbb{R}^3}e^{-\nu(v)(t-s)} \Big|(K^mf)(s,x-v(t-s),v)\Big|dvds\nonumber\\
&\quad+\int_0^t\int_{\mathbb{R}^3}e^{-\nu(v)(t-s)} \Big|(K^cf)(s,x-v(t-s),v)\Big|dvds\nonumber\\
&\quad+\int_0^t\int_{\mathbb{R}^3}e^{-\nu(v)(t-s)} \Big|\Gamma(f,f)(s,x-v(t-s),v)\Big|dvds\nonumber\\
&:=
\int_{\mathbb{R}^3}e^{-\nu(v)t}|f_0(x-vt,v)|dv+H_1+H_2+H_3.
\end{align}
For $H_1$, it follows from \eqref{2.31} that 
\begin{equation}\label{5.5-1}
H_1\leq Cm^{3+\g}\sup_{0\leq s\leq t}\|f(s)\|_{L^\infty}\int_0^{t}\int_{\mathbb{R}^3}e^{-\nu(v)(t-s)} e^{-\f{|v|^2}{10}} dvds\leq Cm^{3+\g}\sup_{0\leq s\leq t}\|f(s)\|_{L^\infty}.
\end{equation}
For $H_2$, we    notice that 
\begin{align}\label{5.5}
H_2
&= \int_{t-\l}^t\int_{\mathbb{R}^3}e^{-\nu(v)(t-s)} \Big|\int_{\mathbb{R}^3_{v'}}l(v,v')f(s,x-v(t-s),v')dv'\Big|dvds\nonumber\\
&\quad+\int_0^{t-\l}\int_{\mathbb{R}^3}e^{-\nu(v)(t-s)} \Big|\int_{\mathbb{R}^3_{v'}}l(v,v')f(s,x-v(t-s),v')dv'\Big|dvds
:=
H_{21}+H_{22}.
\end{align}
It is straightforward to obtain that for $\b> 2$,
\begin{align}
H_{21}&\leq \sup_{0\leq s\leq t}\|h(s)\|_{L^\infty}\int_{t-\l}^t\int_{\mathbb{R}^3}e^{-\nu(v)(t-s)}(1+|v|)^{-\b} \int_{\mathbb{R}^3_{v'}}\Big|l_{w_\b}(v,v')\Big| dv' dv ds\nonumber\\
&\leq C\l \sup_{0\leq s\leq t}\|h(s)\|_{L^\infty},\notag
\end{align}
where we have used \eqref{4.20-1} in the last inequality. For the term $H_{22}$, one notices that 
\begin{align}
H_{22}&=\int_0^{t-\l}\int_{|v|\geq N}e^{-\nu(v)(t-s)} \Big|\int_{\mathbb{R}^3_{v'}}l(v,v')f(s,x-v(t-s),v')dv'\Big|dvds\nonumber\\
&~~~+\int_0^{t-\l}\int_{|v|\leq N}e^{-\nu(v)(t-s)} \Big|\int_{|v'|\geq 2N}l(v,v')f(s,x-v(t-s),v')dv'\Big|dvds\nonumber\\
&~~~+\int_0^{t-\l}\int_{|v|\leq N}e^{-\nu(v)(t-s)} \Big|\int_{|v'|\leq 2N}l(v,v')f(s,x-v(t-s),v')dv'\Big|dvds\nonumber\\
&:=
H_{221}+H_{222}+H_{223}.\notag
\end{align}
It follows from \eqref{4.20-1} and \eqref{4.31} that for $\b>2$, 
\begin{align}
H_{221}&\leq \sup_{0\leq s\leq t}\|h(s)\|_{L^\infty}\cdot\int_0^{t-\l}\int_{|v|\geq N}e^{-\nu(v)(t-s)}w_\b(v)^{-1} \Big|\int_{\mathbb{R}^3_{v'}}l_{w_\b}(v,v')dv'\Big|dvds\nonumber\\
&\leq C_m\sup_{0\leq s\leq t}\|h(s)\|_{L^\infty}\int_0^{t-\l}\int_{|v|\geq N}e^{-\nu(v)(t-s)}(1+|v|)^{-2-\b}\nu(v) dvds\nonumber\\
&\leq \f{C_m}{N^{\b-1}}\sup_{0\leq s\leq t}\|h(s)\|_{L^\infty}\leq \f{C_m}{N}\sup_{0\leq s\leq t}\|h(s)\|_{L^\infty},\notag
\end{align}
and
\begin{align}
H_{222}&\leq e^{-\f{N^2}{20}}\sup_{0\leq s\leq t}\|h(s)\|_{L^\infty}\cdot\int_0^{t-\l}\int_{|v|\leq N}e^{-\nu(v)(t-s)}w_{\b}(v)^{-1}
\Big|\int_{|v'|\geq2N}l_{w_\b}(v,v')e^{\f{|v-v'|^2}{20}}dv'\Big|dvds\nonumber\\
&\leq C_m e^{-\f{N^2}{20}}\sup_{0\leq s\leq t}\|h(s)\|_{L^\infty}.\notag
\end{align}

Since $l_{w_{\b}}(v,v')$ has a possible singularity of $\f1{|v-v'|}$, as before we  choose the smooth compact support function $l_N(v,v')$ satisfying \eqref{4.35}. 
Then it follows from \eqref{4.39} and \eqref{4.35} that 
\begin{align}\label{5.10}
H_{223}&\leq \int_0^{t-\l}\int_{|v|\leq N}e^{-\nu(v)(t-s)}w_{\b}(v)^{-1} \nonumber\\
&\quad\times\int_{|v'|\leq 2N}|l_{w_{\b}}(v,v')-l_N(v,v')|\cdot|h(s,x-v(t-s),v')|dv'dvds\nonumber\\
&\quad+\int_0^{t-\l}\int_{|v|\leq N}e^{-\nu(v)(t-s)} \int_{|v'|\leq 2N}\Big|l_N(v,v')h(s,x-v(t-s),v')\Big|dv'dvds\nonumber\\
&\leq \f{C_m}{N}\sup_{0\leq s\leq t}\|h(s)\|_{L^\infty}+C_N\int_0^{t-\l}e^{-c_N(t-s)}\int_{|v|\leq N,|v'|\leq 2N} \Big|h(s,x-v(t-s),v')\Big|dv'dvds\nonumber\\
&\leq \f{C_m}{N}\sup_{0\leq s\leq t}\|h(s)\|_{L^\infty}+C_N\l^{-\f32}\sqrt{\mathcal{E}(F_0)}+C_N\l^{-3}\mathcal{E}(F_0).
\end{align}
Hence, combining \eqref{5.5}-\eqref{5.10}, one obtains that for $\b>2$, 
\begin{align}\label{5.11}
H_2\leq C_m\left(\l+\f{1}{N}\right)\sup_{0\leq s\leq t}\|h(s)\|_{L^\infty}+C_N\l^{-\f32}\sqrt{\mathcal{E}(F_0)}+C_N\l^{-3}\mathcal{E}(F_0).
\end{align}

Next we estimate $H_3$. Firstly, we note that 
\begin{align}\label{5.12}
H_3&\leq \int_0^t\int_{\mathbb{R}^3}e^{-\nu(v)(t-s)} \int_{\mathbb{R}^3}\int_{\mathbb{S}^2}B(v-u,\t)e^{-\f{|u|^2}{4}}\nonumber\\
&\qquad\qquad\times\Big|f(s,x-v(t-s),u)f(s,x-v(t-s),v)\Big|dud\omega dvds\nonumber\\
&\quad+\int_0^t\int_{\mathbb{R}^3}e^{-\nu(v)(t-s)} \int_{\mathbb{R}^3}\int_{\mathbb{S}^2}B(v-u,\t)e^{-\f{|u|^2}{4}}\nonumber\\
&\qquad\qquad\times\Big|f(s,x-v(t-s),u')f(s,x-v(t-s),v')\Big|dud\omega dvds
:=
H_{31}+H_{32}.
\end{align}
For $H_{31}$, one has that for $\b>\max\{3,3+\g\}$,  
\begin{align}\label{5.13}
H_{31}&\leq C\int_{t-\l}^t\int_{\mathbb{R}^3}e^{-\nu(v)(t-s)}\nu(v)w_{\b}(v)^{-1}\|h(s)\|^2_{L^\infty} dvds\nonumber\\
&\quad+C\int_0^{t-\l}\int_{\mathbb{R}^3}\int_{\mathbb{R}^3}\|h(s)\|_{L^\infty}e^{-\nu(v)(t-s)} w_{\b}(v)^{-1}|v-u|^{\g}e^{-\f{|u|^2}{4}}|f(s,x-v(t-s),u)|dudvds\nonumber\\
&\leq C\l\sup_{0\leq s\leq t}\|h(s)\|^2_{L^\infty}+C\int_0^{t-\l}\Big\{\int_{|v|\geq N}\int_{\mathbb{R}^3_{u}}+\int_{\mathbb{R}^3_v}\int_{|u|\geq N}\Big\}\{\cdots\}dudvds\nonumber\\
&\quad+C\int_0^{t-\l}\int_{|v|\leq N}\int_{|u|\leq N}\{\cdots\}dudvds\nonumber\\
&\leq C\left(\l+\f{1}{N^{\b-3}}\right)\sup_{0\leq s\leq t}\|h(s)\|^2_{L^\infty}+C\int_0^{t-\l}\int_{|v|\leq N}\int_{|u|\leq N}\{\cdots\}dudvds.
\end{align}
To estimate the last term on the RHS of above, for $p>1$  defined in \eqref{4.10},  it follows from the H\"older inequality that
\begin{align}\label{5.14}
&C\int_0^{t-\l}\int_{|v|\leq N}\int_{|u|\leq N}\{\cdots\}dudvds\nonumber\\
&\leq C\int_0^{t-\l}e^{-c_N(t-s)}\|h(s)\|_{L^\infty}\int_{|v|\leq N}\int_{|u|\leq N} w_{\b}(v)^{-1}|v-u|^{\g}e^{-\f{|u|^2}{4}}|f(s,x-v(t-s),u)|dudvds\nonumber\\
&\leq C_N\int_0^{t-\l}e^{-c_N(t-s)}\|h(s)\|_{L^\infty}
\left(\int_{|v|\leq N}\int_{|u|\leq N}|f(s,x-v(t-s),u)|^{\f{p}{p-1}}{du dv}\right)^{1-\f1p}ds\nonumber\\
&\leq C_N\l^{-3}\Big(\sqrt{\mathcal{E}(F_0)}+\mathcal{E}(F_0)\Big)^{1-\f1p}\sup_{0\leq s\leq t}\|h(s)\|^{1+\f1p}_{L^\infty},
\end{align}
where in the last inequality we have used the following fact  that
\begin{align}\label{5.15}
&\int_{|v|\leq N,|u|\leq 3N} |f(s,x-v(t-s),u)|dudv
\leq C_N\f{1+(t-s)^{\f32}}{(t-s)^{\f32}}\sqrt{\mathcal{E}(F_0)}+C_N\f{1+(t-s)^{3}}{(t-s)^{3}}\mathcal{E}(F_0).
\end{align}
Hence, from \eqref{5.13}-\eqref{5.14}, one obtains, for $\b>\max\{3,3+\g\}$,  that 
\begin{equation}\label{5.16}
H_{31}\leq C\left(\l+\f{1}{N^{\b-3}}\right)\sup_{0\leq s\leq t}\|h(s)\|^2_{L^\infty}+C_N\l^{-3}\Big(\sqrt{\mathcal{E}(F_0)}+\mathcal{E}(F_0)\Big)^{1-\f1p}\sup_{0\leq s\leq t}\|h(s)\|^{1+\f1p}_{L^\infty}.
\end{equation}
For $H_{32}$, we notice, for  $\b>\max\{3,3+\g\}$,  that 
\begin{align}
H_{32}&\leq C\int_{t-\l}^t\int_{\mathbb{R}^3}\int_{\mathbb{R}^3}e^{-\nu(v)(t-s)}|v-u|^{\g}w_{\b}(v)^{-1} e^{-\f{|u|^2}{4}}dudvds\cdot\sup_{0\leq s\leq t}\|h(s)\|^2_{L^\infty}\nonumber\\
&~~+\int_0^{t-\l}\int_{\mathbb{R}^3}\int_{\mathbb{R}^3}\int_{\mathbb{S}^2}B(v-u,\t)e^{-\nu(v)(t-s)} w_{\b}(v)^{-1}\nonumber\\
&\qquad\qquad\qquad\times e^{-\f{|u|^2}{4}}\|h(s)\|_{L^\infty}|h(s,x-v(t-s),v')|dud\omega dvds\nonumber\\
&\leq C\l\sup_{0\leq s\leq t}\|h(s)\|^2_{L^\infty}+\int_0^{t-\l}\Big\{\int_{|v|\geq N}\int_{\mathbb{R}^3_{u}}\int_{\mathbb{S}^2}+\int_{\mathbb{R}^3_v}\int_{|u|\geq N}\int_{\mathbb{S}^2}\Big\}\{\cdots\}dud\omega dvds\nonumber\\
&\quad+\int_0^{t-\l}\int_{|v|\leq N}\int_{|u|\leq N}\int_{\mathbb{S}^2}\{\cdots\}dud\omega dvds\nonumber\\
&\leq C\left(\l+\f{1}{N^{\b-3}}\right)\sup_{0\leq s\leq t}\|h(s)\|^2_{L^\infty}+\int_0^{t-\l}\int_{|v|\leq N}\int_{|u|\leq N}\int_{\mathbb{S}^2}\{\cdots\}dud\omega dvds.\notag
\end{align}
To estimate the last term on the RHS of above, we  utilize the changing of variables \eqref{3.9}, \eqref{3.10} and  \eqref{3.11} to obtain that 
\begin{align}\label{5.19}
&\int_0^{t-\l}\int_{|v|\leq N}\int_{|u|\leq N}\int_{\mathbb{S}^2}\{\cdots\}dud\omega dvds\nonumber\\
&\leq C\int_0^{t-\l}e^{-c_N(t-s)}\|h(s)\|_{L^\infty}\int_{|v|\leq N}\int_{|z|\leq 2N}\int_{\mathbb{S}^2}w_{\b}(v)^{-1}\f{|z_{\shortparallel}|}{(|z_{\perp}|+|z_{\shortparallel}|)^{1-\g}}e^{-\f{|v+z_{\shortparallel}+z_{\perp}|^2}{4}}\nonumber\\
&\qquad\qquad\times \Big|h(s,x-v(t-s),v+z_{\shortparallel})\Big|dzd\omega dvds\nonumber\\
&\leq C\int_0^{t-\l}e^{-c_N(t-s)}\|h(s)\|_{L^\infty}{\int_{|v|\leq N}\int_{|\eta|\leq 3N}\int_{z_{\perp}}w_{\b}(v)^{-1}\f{|z_{\perp}|^{\f{\g-1}{2}}}{|\eta-v|^{\f{3-\g}{2}}}e^{-\f{|\eta+z_{\perp}|^2}{4}}}\nonumber\\
&\qquad\qquad{\times \Big|h(s,x-v(t-s),\eta)\Big|dz_{\perp}d\eta dvds}\nonumber\\
&\leq C_N \int_0^{t-\l}e^{-c_N(t-s)}\|h(s)\|^{1+\f1p}_{L^\infty}\left(\int_{|v|\leq N}\int_{|\eta|\leq 3N}|f(s,x-v(t-s),\eta)|d\eta dv\right)^{1-\f1p}ds\nonumber\\
&\leq C_N\l^{-3}\Big(\sqrt{\mathcal{E}(F_0)}+\mathcal{E}(F_0)\Big)^{1-\f1p}\cdot\sup_{0\leq s\leq t}\|h(s)\|^{1+\f1p}_{L^\infty},
\end{align}
where we have used \eqref{5.15} in the last inequality and recall that {$p>1$ is defined in \eqref{4.10}.} 
Hence, combining \eqref{5.12} and \eqref{5.16}-\eqref{5.19},   one obtains that  for $\b>\max\{3,3+\g\}$, 
\begin{equation}\label{5.20}
H_{3}\leq C\Big(\l+\f{1}{N^{\b-3}}\Big)\sup_{0\leq s\leq t}\|h(s)\|^2_{L^\infty}+C_N\l^{-3}\Big(\sqrt{\mathcal{E}(F_0)}+\mathcal{E}(F_0)\Big)^{1-\f1p}\cdot\sup_{0\leq s\leq t}\|h(s)\|^{1+\f1p}_{L^\infty}.
\end{equation}
Submitting \eqref{5.5-1}, \eqref{5.11} and \eqref{5.20}  into \eqref{5.3}, one proves \eqref{5.2}  for $\b> \max\{3, 3+\g\}$. Hence the proof of Lemma \ref{lem5.1} is completed. $\hfill\Box$


\subsection{Global existence and uniqueness}

Now we are in a position to give the 

\medskip

\noindent{\bf Proof of Theorem \ref{thm1.1}:}  Let $\b>\max\{3, 3+\g\}$. In terms of \eqref{4.5}, we make the {\it a priori} assumption
\begin{align}\label{5.22} 
\|h(t)\|_{L^\infty}\leq 2A_0:=
2C_1\Big\{2\bar{M}^2+\sqrt{\mathcal{E}(F_0)}+\mathcal{E}(F_0)\Big\},
\end{align}
where the positive constant $C_1\geq1$ is defined in Lemma \ref{lem4.1}. Then it follows from Lemma \ref{lem4.1} and the a priori assumption \eqref{5.22} that 
\begin{align}\label{5.23}
\|h(t)\|_{L^\infty}\leq A_0+C_1(2A_0)^{\f{9p+1}{5p}}\cdot\sup_{t_1\leq s\leq t,~y\in\Omega} \Big(\int_{\mathbb{R}^3}|f(s,y,\eta)|d\eta\Big)^{\f{p-1}{5p}}.
\end{align}
To estimate the second term on the RHS of \eqref{5.23}, we first notice that for $\Omega=\mathbb{R}^3$ and $t\geq t_1$,
\begin{align}\label{5.40}
\int_{\mathbb{R}^3}e^{-\nu(v)t}|f_0(x-vt,v)|dv&
\leq t_1^{-3}\|f_0\|_{L^1_xL^\infty_v}\leq C\bar{M}^3\|f_0\|_{L^1_xL^\infty_v}.
\end{align}
For $\Omega=\mathbb{T}^3$ and $t\geq t_1$, it holds  that 
\begin{align}\label{5.41}
\int_{\mathbb{R}^3}e^{-\nu(v)t}|f_0(x-vt,v)|dv
&\leq \left(\int_{|v|\geq \tilde{N}}+\int_{|v|\leq \tilde{N}}\right)|f_0(x-vt,v)|dv\notag\\
&\leq  C \|w_\beta f_0\|_{L^\infty}^{\frac{3}{\beta}}\|f_0\|_{L^1_xL^\infty_v}^{1-\frac{3}{\beta}} +\frac{C}{t_1^3}\|f_0\|_{L^1_xL^\infty_v}\notag\\
&\leq C{\bar{M}}^{\frac{3}{\beta}}\|f_0\|_{L^1_xL^\infty_v}^{1-\frac{3}{\beta}} +C \bar{M}^3\|f_0\|_{L^1_xL^\infty_v},
\end{align}
where we have chosen $\tilde{N}=\|w_\beta f_0\|_{L^\infty}^{\f1\b}\|f_0\|_{L^1_xL^\infty}^{-\f1\b}$. 
Then it follows from  Lemma \ref{lem5.1}, \eqref{5.40}, \eqref{5.41} and the a priori assumption \eqref{5.22} that 
\begin{align}\label{5.24}
&\sup_{t_1\leq s\leq t,y\in\mathbb{R}^3}\int_{\mathbb{R}^3}|f(t,y,v)|dv\nonumber\\
&\leq 
\begin{cases}
C\bar{M}^3\|f_0\|_{L^1_xL^\infty_v},~\mbox{for}~\Omega=\mathbb{R}^3\\
C\bar{M}^{\f3\b}\|f_0\|_{L^1_xL^\infty}^{1-\f3\b}{+C\bar{M}^3\|f_0\|_{L^1_xL^\infty_v}}~\mbox{for}~\Omega=\mathbb{T}^3
\end{cases}
+C_N\l^{-\f32}\sqrt{\mathcal{E}(F_0)}+C_N\l^{-3}\mathcal{E}(F_0)\nonumber\\
&\quad +C\Big\{m^{3+\g}+C_m[\l+\f{1}{N}+\f{1}{N^{\b-3}}]\Big\}(2A_0)^2+C_N\l^{-3}\Big(\sqrt{\mathcal{E}(F_0)}+\mathcal{E}(F_0)\Big)^{1-\f1p}(2A_0)^{1+\f1p}.
\end{align}
Note $\b>\max\{3, 3+\g\}$ and $p>1$. One can firstly choose $\l$ sufficiently small, then $N\geq1$ large enough, and finally let $\mathcal{E}(F_0)+\|f_0\|_{L^1_xL^\infty_v}\leq \v_0$ with $\v_0$ small depending only on $\b,\g$ and $\bar{M}$, such that 
\begin{align}\label{5.26}
 4C_1A_0^{\f{4p+1}{5p}}\cdot\sup_{t_1\leq s\leq t,~y\in\Omega} \Big(\int_{\mathbb{R}^3}|f(s,y,\eta)|d\eta\Big)^{\f{p-1}{5p}}\leq \f34,
\end{align}
which together with \eqref{5.23}, yield immediately that
\begin{equation}\label{5.27}
\|h(t)\|_{L^\infty}\leq \f74A_0,
\end{equation}
for all $t\geq 0$. 
Hence we have closed the a priori assumption \eqref{5.22}. Therefore the proof of Theorem \ref{thm1.1} is completed. $\hfill\Box$


\subsection{Positive Lower Bound of Density} At the end of this section, we give the proof of Corollary \ref{cor1.3}. Noting $C_1\geq 1$ and $A_0\geq 1$,  it follows from \eqref{5.26} that for $t\geq t_1$,
\begin{equation}\label{5.28}
\int_{\mathbb{R}^3}|f(t,x,v)|dv\leq \f34.
\end{equation}
Then, using \eqref{5.28}, it is straightforward to get that
\begin{align}
|\r(t,x)-1|=\left|\int_{\mathbb{R}^3}[F(t,x,v)-\mu(v)]dv\right|\leq \int_{\mathbb{R}^3}|f(t,x,v)|dv\leq \f34,\notag
\end{align}
which yields immediately that initial vacuum of the density function should disappear for $t\geq T_0:=
t_1$.  Therefore the proof of Corollary \ref{cor1.3} is completed. $\hfill\Box$

\section{Time-Decay Estimates in Torus}

In this section, we consider the time-decay estimates for the global  solutions  obtained in Theorem \ref{thm1.1}.
Let $\Omega=\mathbb{T}^3$,  we consider the following linearized Boltzmann equation
\begin{align}\label{6.1}
\zeta_t+v\cdot\nabla_x\z+\nu(v)\z-K\z=0,\quad \z(0,x,v)=\z_0(x,v).
\end{align}
Denoting the semigroup of \eqref{6.1} by  $S(t)$, it holds that 
\begin{equation}
\z(t)=S(t)\z_0.\notag
\end{equation}
Let $\z(t,x,v)$ be the solution of the linearized equation \eqref{6.1}, and denote 
\begin{equation}
\xi(t,x,v)=w_\b(v)\z(t,x,v).\notag
\end{equation}
Then it follows from \eqref{6.1} that 
\begin{align}\label{6.1-1}
\xi_{t}+v\cdot\nabla_x\xi+\nu(v) \xi-K_{w_\b}\xi=0,\quad \xi(0,x,v)=\xi_{0}(x,v).
\end{align}
For later use, we denote the semigroup of \eqref{6.1-1} by $U(t)$, and write the solution as 
\begin{equation}
\xi(t)=U(t)\xi_{0}.\notag
\end{equation}

\subsection{Case of Hard Potentials}
In this subsection, we consider the decay estimate for hard potentials on torus.
The following proposition is a starting point for further geting the exponential decay in $L^\infty$ norm.

\begin{proposition}[\cite{Kim2014}]\label{prop6.1}
Let $0\leq \g\leq 1$, $\Omega=\mathbb{T}^3$. Let $\z(t,x,v)$ be any solution to the linearized Boltzmann equation \eqref{6.1} and satisfies the conservations of mass \eqref{1.13}, momentum \eqref{1.13-1} and energy \eqref{1.13-2} with $(M_0,J_0,E_0)=(0,0,0)\in \mathbb{R}\times\mathbb{R}^3\times\mathbb{R}$. Then there exists  positive constants $\sigma>0$ and $C>0$ such that
	\begin{align}
	\notag
	\|S(t)\z_0\|_{L^2}=\|\z(t)\|_{L^2}\leq Ce^{-\sigma t}\|\z_0\|_{L^2},
	\end{align}
	for all $t\geq 0$.
\end{proposition}

Utilizing  Proposition \ref{prop6.1}, we can obtain the following $L^\infty$ decay estimate for the linearized Boltzmann equation.

\begin{lemma}\label{lem6.2}
Let $0\leq \g\leq 1$, $\Omega=\mathbb{T}^3$. Let $\z(t,x,v)$ be any solution to the linear Boltzmann equation \eqref{6.1} and satisfies the conservations of mass \eqref{1.13}, momentum \eqref{1.13-1} and energy \eqref{1.13-2} with $(M_0,J_0,E_0)=(0,0,0)\in \mathbb{R}\times\mathbb{R}^3\times\mathbb{R}$. Then there exists  positive constants $0<\sigma_1\leq \sigma$ and $C>0$ such that
\begin{align}\label{6.6}
\|U(t)\xi_0\|_{L^\infty}=\|\xi(t)\|_{L^\infty}\leq Ce^{-\sigma_1 t}\|w_\b \zeta_0\|_{L^\infty}~~\mbox{for}~~t\geq0.
\end{align}
\end{lemma}

\noindent{\bf Proof.}  Notice that via Lemma 19 in \cite{Guo2}, we only need to prove that there exist $\sigma_2>0$, $T_1>0$ 
and $C_{T_1}$ such that 
\begin{equation}\nonumber
\|\xi(T_1)\|_{L^\infty}\leq e^{-\sigma_2T_1}\|\xi_0\|_{L^\infty}+C_{T_1}\int_0^{T_1}\|\zeta(s)\|_{L^2}ds.
\end{equation}
The rest of the proof is similar to Kim \cite{Kim2014}; see also Guo \cite{Guo2}.  Indeed, our case is simpler than \cite{Kim2014,Guo2} since the the characteristic lines in case without forcing are straight lines. Here we omit the details for brevity of presentations.  $\hfill\Box$

\

Based on the above preparations, we utilize  Lemma \ref{lem6.2} to prove Theorem \ref{thm1.2}.\\[1.5mm]
\noindent{\bf Proof of Theorem \ref{thm1.2}:}  Using the semigroup $U(t)$ for the weighted linearized  Boltzmann equation \eqref{6.1-1}, by the Duhamel Principle,   we have the solution formula  for the nonlinear weighted Boltzmann equation \eqref{4.1} as
\begin{align}
h(t)=U(t)h_0+\int_0^tU(t-s)\Big\{w_{\b}\Gamma(f,f)(s)\Big\}ds.\notag
\end{align}
Then it follows from \eqref{6.6} that 
\begin{align}\label{6.9}
\|h(t)\|_{L^\infty}
&\leq Ce^{-\sigma_1 t}\|h_0\|_{L^\infty}+\Big\|\int_0^tU(t-s)\Big\{w_{\b}\Gamma(f,f)(s)\Big\}ds\Big\|_{L^\infty}.
\end{align}
To bound the last term on the RHS of \eqref{6.9}, we notice that 
\begin{align}\label{6.10}
\int_0^tU(t-s)\Big\{w_{\b}\Gamma(f,f)(s)\Big\}ds=&\int_0^te^{-\nu(v)(t-s)}\Big\{w_{\b}\Gamma(f,f)(s)\Big\}ds\nonumber\\
&+\int_0^t\int_s^te^{-\nu(v)(t-s_1)}K_{w_{\b}}\Big\{U(s_1-s)w_{\b}\Gamma(f,f)(s)\Big\}ds_1ds.
\end{align}
For the first term on the RHS of \eqref{6.10}, it follows from \eqref{4.8-1} that 
\begin{align}\label{6.11}
&\Big|\int_0^te^{-\nu(v)(t-s)}\Big\{w_{\b}(v)\Gamma(f,f)(s)\Big\}ds\Big|\nonumber\\
&\leq C\int_0^te^{-\nu(v)(t-s)}\nu(v)\|h(s)\|_{L^\infty}^{\f{9p+1}{5p}}\sup_{y\in\Omega}\Big(\int_{\mathbb{R}^3}|f(s,y,\eta)|d\eta\Big)^{\f{p-1}{5p}}ds\nonumber\\
&\leq C\int_0^te^{-\nu(v)(t-s)}\nu(v) e^{-\f{\sigma_1}2s}\nonumber\\
&\quad\times\sup_{0\leq s\leq t, y\in\Omega}
\left\{\Big[e^{\f{\s_1}2s}\|h(s)\|_{L^\infty}\Big]\cdot\|h(s)\|^{\f{4p+1}{5p}}_{L^\infty}\Big(\int_{\mathbb{R}^3}|f(s,y,\eta)|d\eta\Big)^{\f{p-1}{5p}}\right\}ds\nonumber\\
&\leq Ce^{-\f{\sigma_1}{2} t}\sup_{0\leq s\leq t, y\in\Omega}
\left\{\Big[e^{\f{\s_1}2s}\|h(s)\|_{L^\infty}\Big]\cdot\|h(s)\|^{\f{4p+1}{5p}}_{L^\infty}\Big(\int_{\mathbb{R}^3}|f(s,y,\eta)|d\eta\Big)^{\f{p-1}{5p}}\right\}.
\end{align}

To estimate the second term on the RHS of \eqref{6.10},  as in \cite{Guo2} we define a new semigroup $\tilde{U}(t)$ such that it solves 
\begin{align}
\Big\{\partial_t+v\cdot\nabla_x+\nu(v)-K_{\tilde{w}}\Big\}\{\tilde{U}(t)\tilde{h}_0\}=0,~~\tilde{U}(0)\tilde{h}_0=\tilde{h}_0,\notag
\end{align}
with $\tilde{w}(v)=\f{w_{\b}(v)}{\sqrt{1+|v|^2}}$. A direct calculation shows that 
\begin{equation}\nonumber
\sqrt{1+|v|^2}\tilde{U}(t),
\end{equation}
also solves the original  weighted linearized Boltzmann equation \eqref{6.1-1}. Then the uniqueness in $L^\infty$ class with $\tilde{h}_0=\f{h_0}{\sqrt{1+|v|^2}}$ yields that 
\begin{align}
U(t)h_0\equiv \sqrt{1+|v|^2}\tilde{U}(t)\Big\{\f{h_0}{\sqrt{1+|v|^2}}\Big\}.\notag
\end{align}
Here we  point out that \eqref{6.6} also holds for semigroup $\tilde{U}(t)$.
Then it follows from \eqref{6.6} and \eqref{4.8-1} that 
\begin{align}\label{6.14}
&\Big|\int_0^t\int_s^te^{-\nu(v)(t-s_1)}K_{w_{\b}}\Big\{U(s_1-s)w_{\b}\Gamma(f,f)(s)\Big\}ds_1ds\Big|\nonumber\\
&\leq \int_0^t\int_s^te^{-\nu_0(t-s_1)}\Big|\int_{\mathbb{R}^3_{v'}}k_{w_{\b}}(v,v')\sqrt{1+|v'|^2}\Big\{\tilde{U}(s_1-s)\f{w_{\b}}{\sqrt{1+|v|^2}}\Gamma(f,f)(s)\Big\}dv'\Big|ds_1ds\nonumber\\
&\leq \int_0^t\int_s^te^{-\nu_0(t-s_1)}\Big|\int_{\mathbb{R}^3_{v'}}k_{w_{\b}}(v,v')\sqrt{1+|v'|^2}dv'\Big|\nonumber\\
&\qquad\qquad\qquad\qquad\qquad \times\Big\|\Big\{\tilde{U}(s_1-s)\f{w_{\b}}{\sqrt{1+|v|^2}}\Gamma(f,f)(s)\Big\}\Big\|_{L^\infty}ds_1ds\nonumber\\
&\leq \int_0^t\int_s^te^{-\nu_0(t-s_1)}e^{-\sigma_1(s_1-s)}\Big\|\f{w_{\b}}{\sqrt{1+|v|^2}}\Gamma(f,f)(s)\Big\|_{L^\infty}ds_1 ds\nonumber\\
&\leq C\int_0^t\int_s^te^{-\nu_0(t-s_1)}e^{-\sigma_1(s_1-s)}\Big\|\f{\nu(v)}{\sqrt{1+|v|^2}}\Big\|_{L^\infty}\|h(s)\|_{L^\infty}^{\f{9p+1}{5p}}\sup_{y\in\Omega}\Big(\int_{\mathbb{R}^3}|f(s,y,\eta)|d\eta\Big)^{\f{p-1}{5p}}ds_1ds\nonumber\\
&\leq Ce^{-\f{\sigma_1}{2} t}\sup_{0\leq s\leq t, y\in\Omega}
\left\{\Big[e^{\f{\s_1}2s}\|h(s)\|_{L^\infty}\Big]\cdot\|h(s)\|^{\f{4p+1}{5p}}_{L^\infty}\Big(\int_{\mathbb{R}^3}|f(s,y,\eta)|d\eta\Big)^{\f{p-1}{5p}}\right\}.
\end{align}
Combining \eqref{6.9}-\eqref{6.11} and \eqref{6.14} and using \eqref{5.27}, one obtains that 
\begin{align}\label{6.15}
&\sup_{0\leq s\leq t}\Big\{e^{\f{\s_1}{2}s}\|h(s)\|_{L^\infty}\Big\}\nonumber\\
&\leq C\|h_0\|_{L^\infty}+C\sup_{0\leq s\leq t, y\in\Omega}
\left\{\Big[e^{\f{\s_1}2s}\|h(s)\|_{L^\infty}\Big]\cdot\|h(s)\|^{\f{4p+1}{5p}}_{L^\infty}\Big(\int_{\mathbb{R}^3}|f(s,y,\eta)|d\eta\Big)^{\f{p-1}{5p}}\right\}\nonumber\\
&\leq C\Big\{\|h_0\|_{L^\infty}+\sup_{0\leq s\leq 1}\|h(s)\|^2_{L^\infty}\Big\}\nonumber\\
&\quad+C\sup_{1\leq s\leq t}\Big[e^{\f{\s_1}2s}\|h(s)\|_{L^\infty}\Big]
\cdot\sup_{1\leq s\leq t,y\in\Omega}\left\{\|h(s)\|^{\f{4p+1}{5p}}_{L^\infty}\Big(\int_{\mathbb{R}^3}|f(s,y,\eta)|d\eta\Big)^{\f{p-1}{5p}}\right\}\nonumber\\
&\leq C_2\bar{M}^4+C_2\sup_{1\leq s\leq t}\Big[e^{\f{\s_1}2s}\|h(s)\|_{L^\infty}\Big]
\cdot\sup_{1\leq s\leq t,y\in\Omega}\left\{\|h(s)\|^{\f{4p+1}{5p}}_{L^\infty}\Big(\int_{\mathbb{R}^3}|f(s,y,\eta)|d\eta\Big)^{\f{p-1}{5p}}\right\}.
\end{align}
Then, using \eqref{5.24} and  similar arguments as in \eqref{5.26},  if $\v_0$ is small enough, one can obtain that 
\begin{align}\label{6.20}
&C_2\sup_{1\leq s\leq t,y\in\Omega}\left\{\|h(s)\|^{\f{4p+1}{5p}}_{L^\infty}\Big(\int_{\mathbb{R}^3}|f(s,y,\eta)|d\eta\Big)^{\f{p-1}{5p}}\right\}\leq \f12.
\end{align}
Substituting \eqref{6.20} into \eqref{6.15}, one gets that 
\begin{align}\label{6.21}
e^{\f{\sigma_1}{2}t}\|h(t)\|_{L^\infty}
\leq 2C_2\bar{M}^4,\quad \forall\,t\geq0.
\end{align}
Finally,   choosing
$$\sigma_0=\f{\sigma_1}{2}~~~\mbox{and}~~~\tilde{C}_2=2C_2\bar{M}^4,$$
we then obtain \eqref{1.29} from \eqref{6.21}.  Therefore the proof of Theorem \ref{thm1.2} is completed. $\hfill\Box$

\subsection{Case of Soft Potentials}

In this subsection, we consider the decay estimates for soft potentials on torus.  Firstly, we define the Fourier transformation as
\begin{align}\notag
\hat{\zeta}(t,k,v)=\int_{\Omega}e^{-ik\cdot x} \zeta(t,x,v)dx,\quad k\in\mathbb{Z}^3.
\end{align}
Then the have the following estimate, whose proof can be found in \cite{Duan-Yang-Zhao,Strain}.

\begin{proposition}[\cite{Duan-Yang-Zhao}]\label{prop8.1}
Let $-3< \g<0$, and let $d\geq0$, $r>0$ be given constants.  Let $\zeta(t,x,v)$ be any solution to the linearized Boltzmann equation \eqref{6.1} and satisfies the conservations of mass \eqref{1.13}, momentum \eqref{1.13-1} and energy \eqref{1.13-2} with $(M_0,J_0,E_0)=(0,0,0)\in \mathbb{R}\times\mathbb{R}^3\times\mathbb{R}$. 
Then the following estimate holds
	\begin{align}
	\|\sqrt{\nu(\cdot)}^{-d}\hat{\zeta}(t,k,\cdot)\|^2_{L^2_{v}}\leq C(1+t)^{-r}\|\sqrt{\nu(\cdot)}^{-d-r_+}\hat{\zeta}_0(k,\cdot)\|^2_{L^2_v},\notag
	\end{align}
for all  $t\geq0$ and $k\in\mathbb{Z}^3$,	where $r_+$ denotes the arbitrary constant which is strictly greater than $r$.
\end{proposition}

Using Proposition \ref{prop8.1} and Plancherel theorem, we have  the following $L^2$ decay estimate.

\begin{lemma}
Under the assumptions of Proposition \ref{prop8.1}, the  following estimate holds
\begin{align}\label{8.2}
\|\sqrt{\nu}^{-d}S(t)\zeta_0\|^2_{L^2}=\|\sqrt{\nu}^{-d}\zeta(t)\|^2_{L^2}\leq C(1+t)^{-r}\|\sqrt{\nu}^{-d-r_+}\zeta_0\|^2_{L^2},
\end{align}
for all  $t\geq0$.
\end{lemma}


We have the following $L^\infty$-decay estimate for the solutions to the linearized Boltzmann equation.
\begin{lemma}\label{lem8.3}
Under the assumptions of  Proposition \ref{prop8.1},  it holds  that 
\begin{align}\label{8.5}
\|S(t)\z_0\|_{L^\infty}=\|\z(t)\|_{L^\infty}\leq C(1+t)^{-r}\|w_{2+|\g|r}\z_0\|_{L^\infty},
\end{align}
for any given $r\in(0,1+\f{2}{|\g|})$.
\end{lemma}

\noindent{\bf Proof.} It is noted that 
\begin{align}\label{8.6-1}
\z(t,x,v)&= e^{-\nu(v)t}\z_0(x-vt,v)+\int_0^t e^{-\nu(v)(t-s)}K^m\z(s,x-v(t-s),v)ds\nonumber\\
&\quad +\int_0^t e^{-\nu(v)(t-s)}K^c\z(s,x-v(t-s),v)ds,
\end{align}
which yields immediately that 
\begin{align}
|\z(t,x,v)|&\leq e^{-\nu(v)t}|\z_0(x-vt,v)|+\int_0^t e^{-\nu(v)(t-s)}|K^m\z(s,x-v(t-s),v)|ds\nonumber\\
&\quad +\int_0^t e^{-\nu(v)(t-s)}|K^c\z(s,x-v(t-s),v)|ds\nonumber\\
&:=
L_1+L_2+L_3,\notag
\end{align} 
where $m>0$ is a small constant to be chosen later.  

Firstly, it is easy to get that 
\begin{align}\label{8.8}
L_1\leq C(1+t)^{-r}\|\nu^{-r}\z_0\|_{L^\infty}.
\end{align}
It follows from \eqref{2.31} that
\begin{align}
L_2&\leq Cm^{3+\g}\int_0^t e^{-\nu(v)(t-s)}e^{-\f{|v|^2}{10}}\|\z(s)\|_{L^\infty}ds\nonumber\\
&\leq Cm^{3+\g}e^{-\f{|v|^2}{20}}\cdot\int_0^t(1+t-s)^{-1-r}\|\z(s)\|_{L^\infty}ds\nonumber\\
&\leq Cm^{3+\g}e^{-\f{|v|^2}{20}}\sup_{0\leq s\leq t}\Big\{(1+s)^r\|\z(s)\|_{L^\infty}\Big\}\cdot\int_0^t(1+t-s)^{-1-r}(1+s)^{-r}ds\nonumber\\
&\leq Cm^{3+\g}e^{-\f{|v|^2}{20}}(1+t)^{-r}\sup_{0\leq s\leq t}\Big\{(1+s)^r\|\z(s)\|_{L^\infty}\Big\}.\label{8.7}
\end{align}
To bound  $L_3$, we use \eqref{8.6-1} again to get that 
\begin{align}\label{8.9}
L_3&\leq \int_0^t e^{-\nu(v)(t-s)}\int_{\mathbb{R}^3_{v'}}|l(v,v')|e^{-\nu(v')s}|\z_0(x-v(t-s)-v's,v')|ds\nonumber\\
&\quad+\int_0^t e^{-\nu(v)(t-s)}\int_{\mathbb{R}^3_{v'}}|l(v,v')|\Big\{\int_0^se^{-\nu(v')(s-\tau)}\nonumber\\
&\qquad\qquad\qquad\times|K^m\z(\tau,x-v(t-s)-v'(s-\tau),v')|d\tau\Big\}dv'ds\nonumber\\
&\quad+\int_0^t e^{-\nu(v)(t-s)}\int_{\mathbb{R}^3_{v'}}\int_{\mathbb{R}^3_{v''}}|l(v,v')l(v',v'')|\int_0^se^{-\nu(v')(s-\tau)}\nonumber\\
&\qquad\qquad\qquad\times|\z(\tau,x-v(t-s)-v'(s-\tau),v'')|d\tau dv''dv'ds\nonumber\\
&:=
L_{31}+L_{32}+L_{33}.
\end{align}
For $L_{31}$, it follows from \eqref{2.40} that 
\begin{align}
L_{31}&\leq C\|\nu^{-r}\z_0\|_{L^\infty}\int_0^t e^{-\nu(v)(t-s)}(1+s)^{-r}\cdot\int_{\mathbb{R}_{v'}}|l(v,v')|dv'ds\nonumber\\
&\leq C_m\|\nu^{-r}\z_0\|_{L^\infty}\int_0^t e^{-\nu(v)(t-s)}(1+s)^{-r}\cdot\nu(v)^{1+\f2{|\g|}}ds\nonumber\\
&\leq C_m\|\nu^{-r}\z_0\|_{L^\infty}\int_0^t (1+t-s)^{-1-\f2{|\g|}}(1+s)^{-r}ds\nonumber\\
&\leq C_m(1+t)^{-r}\|\nu^{-r}\z_0\|_{L^\infty}.\notag
\end{align}
Using \eqref{2.31} and \eqref{2.40-1},  one can obtain that 
\begin{align}\label{8.11}
L_{32}&\leq Cm^{3+\g}\int_0^t e^{-\nu(v)(t-s)}\int_{\mathbb{R}^3_{v'}}|l(v,v')|\int_0^se^{-\nu(v')(s-\tau)}e^{-\f{|v'|^2}{10}}\|\z(\tau)\|_{L^\infty}d\tau dv'ds\nonumber\\
&\leq Cm^{3+\g}\sup_{0\leq \tau\leq t}\Big\{(1+\tau)^r\|\z(\tau)\|_{L^\infty}\Big\}\int_0^t (1+t-s)^{-1-r}\nonumber\\
&\quad\times\int_0^s(1+s-\tau)^{-1-r}(1+\tau)^{-r}d\tau ds\nonumber\\
&\leq Cm^{3+\g}(1+t)^{-r}\sup_{0\leq \tau\leq t}\Big\{(1+\tau)^r\|\z(\tau)\|_{L^\infty}\Big\}.
\end{align}
Now we concentrate on the term $L_{33}$. As before, we divide it into several cases.

\medskip
\noindent{\it Case 1.} For $|v|\geq N$, then it follows from \eqref{2.40} that 
\begin{align}\label{8.12-1}
L_{33}
&\leq C_m\sup_{0\leq \tau\leq t}\Big\{(1+\tau)^r\|\z(\tau)\|_{L^\infty}\Big\} \int_0^t e^{-\nu(v)(t-s)} \f{\nu(v)}{(1+|v|)^2}\nonumber\\
&\qquad\times\int_0^s(1+s-\tau)^{-1-\f{2}{|\g|}} (1+\tau)^{-r}d\tau ds\nonumber\\
&\leq \f{C_m}{N^{\d|\g|}}\sup_{0\leq \tau\leq t}\Big\{(1+\tau)^r\|\z(\tau)\|_{L^\infty}\Big\} \int_0^t (1+t-s)^{-1-\f{2}{|\g|}+\d}  \nonumber\\
&\qquad\times\int_0^s(1+s-\tau)^{-1-\f{2}{|\g|}} (1+\tau)^{-r}d\tau ds,
\end{align}
where $\d>0$ is a small positive constant such that 
\begin{equation}\label{8.12-2}
0<r\leq 1+\f{2}{|\g|}-\d~~\mbox{and}~~1+\f{2}{\g}-\d>1.
\end{equation}
It is noted that such $\d>0$ must exist since $r<1+\f{2}{|\g|}$. Then, from \eqref{8.12-2} and  a direct calculation, one can get that 
\begin{align}
\int_0^t (1+t-s)^{-1-\f{2}{|\g|}+\d}  \int_0^s(1+s-\tau)^{-1-\f{2}{|\g|}} (1+\tau)^{-r}d\tau ds\leq C(1+t)^{-r},\notag
\end{align}
which together with \eqref{8.12-1}, yield that 
\begin{align}\label{8.12-4}
L_{33}\leq \f{C_m}{N^{\d|\g|}}(1+t)^{-r}\cdot\sup_{0\leq \tau\leq t}\Big\{(1+\tau)^r\|\z(\tau)\|_{L^\infty}\Big\}. 
\end{align}

\noindent{\it Case 2.}  For either $|v|\leq N,~|v'|\geq2N$ or $|v'|\leq2N,~|v''|\geq 3N$. It is noted that we have either $|v-v'|\geq N$ or $|v'-v''|\geq N$. Then it follows from \eqref{4.31} that
\begin{align}\label{8.12-7}
L_{33}&\leq e^{-\f{|N|^2}{20}}\sup_{0\leq \tau\leq t}\Big\{(1+\tau)^r\|\z(\tau)\|_{L^\infty}\Big\} \int_0^t e^{-\nu(v)(t-s)}\int_{\mathbb{R}^3_{v'}}\Big|l(v,v')e^{\f{|v-v'|^2}{20}}\Big|dv'\nonumber\\
&~~~~~~~~~~~~\times
\int_0^se^{-\nu(v')(s-\tau)}(1+\tau)^{-r}\int_{\mathbb{R}^3_{v''}}\Big|l(v',v'')e^{\f{|v'-v''|^2}{20}}\Big|dv''d\tau ds\nonumber\\
& \leq C_me^{-\f{N^2}{20}}(1+t)^{-r}\sup_{0\leq \tau\leq t}\Big\{(1+\tau)^r\|\z(\tau)\|_{L^\infty}\Big\}.
\end{align}

\

\noindent{\it Case 3.}  $|v|\leq N,~|v'|\leq2N,~|v''|\leq3N$. This is the last remaining case. Firstly, we note that 
\begin{align}\label{8.12-8}
&\int_0^te^{-\nu(v)(t-s)} \int_{|v'|\leq 2N,|v''|\leq3N}|l(v,v')l(v',v'')|\nonumber\\
&\qquad \times \int_0^se^{-\nu(v')(s-\tau)} |\z(\tau,x-v(t-s)-v'(s-\tau),v'')|d\tau dv'dv''ds\nonumber\\
&\leq \int_0^te^{-\nu(v)(t-s)} \int_{|v'|\leq 2N,|v''|\leq3N}|l(v,v')l(v',v'')|\nonumber\\
&\qquad\times \int_{s-\l}^se^{-\nu(v')(s-\tau)} |\z(\tau,x-v(t-s)-v'(s-\tau),v'')|d\tau dv'dv''ds\nonumber\\
&\quad+\int_0^te^{-\nu(v)(t-s)} \int_{|v'|\leq 2N,|v''|\leq3N}|l(v,v')l(v',v'')|\nonumber\\
&\qquad\quad\times \int_0^{s-\l}e^{-\nu(v')(s-\tau)} |\z(\tau,x-v(t-s)-v'(s-\tau),v'')|d\tau dv'dv''ds.
\end{align}
We can bound the first term on the RHS of \eqref{8.12-8} by
\begin{align}
&C_m \l \sup_{0\leq \tau\leq t}\Big\{(1+\tau)^r\|\z(\tau)\|_{L^\infty}\Big\} \int_0^te^{-\nu(v)(t-s)} \nonumber\\
&\qquad\times\int_{|v'|\leq 2N,|v''|\leq3N}|l(v,v')l(v',v'')| (1+s)^{-r}dv'dv''ds\nonumber\\
&\leq C_m \l \sup_{0\leq \tau\leq t}\Big\{(1+\tau)^r\|\z(\tau)\|_{L^\infty}\Big\} \int_0^te^{-\nu(v)(t-s)} \f{\nu(v)^2}{(1+|v|)^4} (1+s)^{-r}ds\nonumber\\
&\leq  C_m\l (1+t)^{-r}\sup_{0\leq \tau\leq t}\Big\{(1+\tau)^r\|\z(\tau)\|_{L^\infty}\Big\}.\notag
\end{align}

Now we shall estimate the second term on the RHS of \eqref{8.12-8}.  Since $l(v,v')$ has singularity of $|v-v'|^{-1}$,  as before, we can choose a smooth compact support function $\tilde{l}_N(v,v')$ such that
\begin{align}\label{8.12-10}
\sup_{|p|\leq 3N}\int_{|v'|\leq 3N}\Big|l(p,v')-\tilde{l}_N(p,v')\Big|dv'\leq C_m N^{2\g-5}.
\end{align}
Noting
\begin{align}\label{8.12-11}
l(v,v')l(v',v'')&=\Big(l(v,v')-\tilde{l}_N(v,v')\Big)l(v',v'')\nonumber\\
&\quad+\Big(l(v',v'')-\tilde{l}_N(v',v'')\Big)\tilde{l}_N(v,v')+\tilde{l}_N(v,v')\tilde{l}_N(v',v''),
\end{align}
and then using \eqref{2.40-1}, \eqref{8.12-10}  and \eqref{8.12-11}, we can bound the second term on the RHS of \eqref{8.12-8} by 
\begin{align}
&C_m N^{2\g-5}\sup_{0\leq \tau\leq t}\Big\{(1+\tau)^{-r}\|\z(\tau)\|_{L^\infty}\Big\}\int_0^te^{-cN^{\g}(t-s)}\int_0^{s-\l}e^{-cN^{\g}(s-\tau)} (1+\tau)^{-r}d\tau ds \nonumber\\
&\quad+\int_0^te^{-cN^{\g}(t-s)} \int_{|v'|\leq 2N,|v''|\leq3N}|\tilde{l}_N(v,v')\tilde{l}_N(v',v'')|\nonumber\\
&\qquad\qquad\times \int_0^{s-\l}e^{-cN^{\g}(s-\tau)} |\z(\tau,x-v(t-s)-v'(s-\tau),v'')|d\tau dv'dv''ds\nonumber\\
&\leq \f{C_m}{N}(1+t)^{-r}\sup_{0\leq s\leq t}\Big\{(1+\tau)^{-r}\|\z(\tau)\|_{L^\infty}\Big\}\nonumber\\
&\quad+C_{N,m}\int_0^t\int_0^{s-\l}e^{-cN^{\g}(t-s)} e^{-cN^{\g}(s-\tau)} \nonumber\\
&\qquad\qquad\qquad\qquad\times\int_{|v'|\leq 2N,|v''|\leq3N} |\z(\tau,x-v(t-s)-v'(s-\tau),v'')|dv'dv''d\tau ds,\notag
\end{align}
where we have used the facts that   $\tilde{l}_N(v,v')\tilde{l}_N(v',v'')$ is bounded and 
\begin{equation}
\nu(v)\geq cN^{\g}~~\mbox{for}~|v|\leq N,~~\mbox{and}~~\nu(v')\geq cN^{\g}~~\mbox{for}~|v'|\leq 2N.\notag
\end{equation}
As in Section 4, using the changing of variables, one can obtains that 
\begin{align}\label{8.12-14} 
& C_{N,m}\int_0^t\int_0^{s-\l}e^{-cN^{\g}(t-s)} e^{-cN^{\g}(s-\tau)} \nonumber\\
&\qquad\qquad\times\int_{|v'|\leq 2N,|v''|\leq3N} |\z(\tau,x-v(t-s)-v'(s-\tau),v'')|dv'dv''d\tau ds
\nonumber\\
&\leq C_{N,m,\l}\int_0^t\int_0^{s-\l}e^{-cN^{\g}(t-s)} e^{-cN^{\g}(s-\tau)} \|\z(\tau)\|_{L^2}d\tau\nonumber\\
&\leq C_{N,m,\l}(1+t)^{-r}\sup_{0\leq\tau\leq t}\Big\{(1+\tau)^r\|\z(\tau)\|_{L^2}\Big\}
\leq C_{N,m}(1+t)^{-r}\|\nu^{-r_+}\z_0\|_{L^2}\nonumber\\
&\leq C_{N,m,\l}(1+t)^{-r}\|w_{2+|\g|r}\z_0\|_{L^\infty},
\end{align} 
where we have used \eqref{8.2} with $d=0$. Thus, combining \eqref{8.12-8}-\eqref{8.12-14}, one gets that 
\begin{align}\label{8.12-15}
&\int_0^te^{-\nu(v)(t-s)} \int_{|v'|\leq 2N,|v''|\leq3N}|l(v,v')l(v',v'')|\nonumber\\
&\qquad\qquad\times \int_0^se^{-\nu(v')(s-\tau)} |\z(\tau,x-v(t-s)-v'(s-\tau),v'')|d\tau dv'dv''ds\nonumber\\
&\leq  C_m\Big(\l +\f{1}{N}\Big)(1+t)^{-r}\cdot\sup_{0\leq \tau\leq t}\Big\{(1+\tau)^r\|\z(\tau)\|_{L^\infty}\Big\}\nonumber\\
&\quad+C_{N,m,\l}(1+t)^{-r}\cdot\|w_{2+|\g|r}\z_0\|_{L^\infty}.
\end{align}

Therefore, it follows from \eqref{8.12-4},  \eqref{8.12-7} and \eqref{8.12-15} that  
\begin{align}
L_{33}&\leq C_m(\l +\f{1}{N}+\f{1}{N^{\d|\g|}})(1+t)^{-r}\cdot\sup_{0\leq \tau\leq t}\Big\{(1+\tau)^r\|\z(\tau)\|_{L^\infty}\Big\}\nonumber\\
&\quad+C_{N,m,\l}(1+t)^{-r}\cdot\|w_{2+|\g|r}\z_0\|_{L^\infty},\notag
\end{align}
where together with \eqref{8.8}, \eqref{8.7} and  \eqref{8.9}-\eqref{8.11}, yield that 
\begin{align}
\sup_{0\leq s\leq t}\Big\{(1+s)^r\|\z(s)\|_{L^\infty}\Big\}&\leq C\Big\{m^{3+\g}+C_m(\l +\f{1}{N}+\f{1}{N^{\d|\g|}})\Big\}\sup_{0\leq \tau\leq t}\Big\{(1+s)^r\|\z(s)\|_{L^\infty}\Big\}\nonumber\\
&\quad+C_{N,m,\l}\|w_{2+|\g|r}\z_0\|_{L^\infty}.\notag
\end{align}
Note $-3<\g<0$. By first choosing $m$ small, then $\l$ small, and finally letting $N$ sufficiently large so that $C\Big\{m^{3+\g}+C_m(\l +\f{1}{N}+\f{1}{N^{\d|\g|}})\Big\}\leq \f12$, one obtains that 
\begin{align}
\notag
\|\z(t)\|_{L^\infty}&\leq C(1+t)^{-r}\|w_{2+|\g|r}\z_0\|_{L^\infty},
\end{align}
for all $t\geq 0$. This yields immediately \eqref{8.5}. Thus we complete the proof of this lemma.    $\hfill\Box$

\medskip

Based on the above preparations, we now use Lemma \ref{lem8.3} to give the

\medskip

\noindent{\bf Proof of Theorem \ref{thm1.3}:}  Using the semigroup $S(t)$ for the linearized  Boltzmann equation \eqref{6.1}, by the Duhamel Principle,   we have the solution formula  for the nonlinear  Boltzmann equation \eqref{1.9} as
\begin{align}
f(t)=S(t)f_0+\int_0^tS(t-s)\Big\{\Gamma(f,f)(s)\Big\}ds.\notag
\end{align}
From now on, we take $r:=
1+\f{2}{|\g|}-\d>1+\f{1}{|\g|}$ with $\d$ being an arbitary small positive constant such that $0<\d\leq \f13$. Then it follows from \eqref{8.5} that 
\begin{multline}\label{8.17}
\|f(t)\|_{L^\infty}
\leq C(1+r)^{-r}\|w_{2+|\g|r}f_0\|_{L^\infty}+C\int_0^t(1+t-s)^{-r}\left\|w_{2+|\g|r}\Big\{\Gamma(f,f)(s)\Big\}\right\|_{L^\infty}ds\\
\leq C(1+r)^{-r}\|w_{4+|\g|}f_0\|_{L^\infty}+C\int_0^t(1+t-s)^{-r}\left\|w_{4+|\g|}\Big\{\Gamma(f,f)(s)\Big\}\right\|_{L^\infty}ds.
\end{multline}
To estimate the last term on the RHS of \eqref{8.17}, we note, from \eqref{4.8-1}, that 
\begin{align}\label{8.20}
&\Big|(1+|v|)^{4+|\g|}\Big\{\Gamma(f,f)(s,x-v(t-s),v)\Big\}\Big|\nonumber\\
&\leq C\|w_4f(s)\|_{L^\infty}\|w_{\f12}f(s)\|^{\f{4p+1}{5p}}_{L^\infty}\sup_{y\in\Omega}\Big(\int_{\mathbb{R}^3}|f(s,y,\eta)|d\eta\Big)^{\f{p-1}{5p}}\nonumber\\
&\leq C\|w_{\b}f(s)\|^{\f4\b}_{L^\infty}\|f(s)\|^{1-\f4\b}_{L^\infty}
\|w_{\b}f(s)\|^{\f1{2\b}\f{4p+1}{5p}}_{L^\infty}\|f(s)\|^{\f{4p+1}{5p}(1-\f1{2\b})}_{L^\infty}\sup_{y\in\Omega}\Big(\int_{\mathbb{R}^3}|f(s,y,\eta)|d\eta\Big)^{\f{p-1}{5p}}\nonumber\\
&\leq C\|f(s)\|^{1-\f4\b+\f{4p+1}{5p}(1-\f1{2\b})}_{L^\infty}\|w_{\b}f(s)\|^{\f4\b+\f{4p+1}{5p}\f1{2\b}}_{L^\infty}\sup_{y\in\Omega}\Big(\int_{\mathbb{R}^3}|f(s,y,\eta)|d\eta\Big)^{\f{p-1}{5p}}\nonumber\\
&\leq C\|f(s)\|_{L^\infty}\|w_{\b}f(s)\|^{\f{4p+1}{5p}}_{L^\infty}\sup_{y\in\Omega3}\Big(\int_{\mathbb{R}^3}|f(s,y,\eta)|d\eta\Big)^{\f{p-1}{5p}},
\end{align} 
where we have used $\f{4p+1}{5p}(1-\f1{4\b})-\f4\b\geq 0$ due to  $\b\geq\f92$ and $1<p\leq \f87$. Then it follows from \eqref{8.20} that 
\begin{align}
&C\int_0^t(1+t-s)^{-r}\left\|(1+|v|)^{4+|\g|}\Big\{\Gamma(f,f)(s)\Big\}\right\|_{L^\infty}ds\nonumber\\
&\leq C\int_0^t(1+t-s)^{-r}(1+s)^{-r}\nonumber\\
&\qquad\times \sup_{0\leq s\leq t,y\in\Omega}\left\{\Big[(1+s)^r\|f(s)\|_{L^\infty}\Big]\|w_{\b}f(s)\|^{\f{4p+1}{5p}}_{L^\infty}\Big(\int_{\mathbb{R}^3}|f(s,y,\eta)|d\eta\Big)^{\f{p-1}{5p}}\right\}ds\nonumber\\
&\leq C(1+t)^{-r}\sup_{0\leq s\leq t,y\in\Omega}\left\{\Big[(1+s)^r\|f(s)\|_{L^\infty}\Big]\|w_{\b}f(s)\|^{\f{4p+1}{5p}}_{L^\infty}\Big(\int_{\mathbb{R}^3}|f(s,y,\eta)|d\eta\Big)^{\f{p-1}{5p}}\right\},\nonumber
\end{align}
which together with \eqref{5.27} and \eqref{8.17}, yield that for $\b\geq\max\{\f92,4+|\g|\}$, 
\begin{align}\label{8.25}
&\sup_{0\leq s\leq t}\Big\{(1+s)^r\|f(s)\|_{L^\infty}\Big\}\nonumber\\
&\leq C\|w_\b f_0\|_{L^\infty}+\sup_{0\leq s\leq t,y\in\Omega}\left\{\Big[(1+s)^r\|f(s)\|_{L^\infty}\Big]\|w_{\b}f(s)\|^{\f{4p+1}{5p}}_{L^\infty}\Big(\int_{\mathbb{R}^3}|f(s,y,\eta)|d\eta\Big)^{\f{p-1}{5p}}\right\}\nonumber\\
&\leq C_3\bar{M}^4+C_3\sup_{1\leq s\leq t}\Big[(1+s)^r\|f(s)\|_{L^\infty}\Big]\sup_{1\leq s\leq t,y\in\Omega}\left\{\|w_{\b}f(s)\|^{\f{4p+1}{5p}}_{L^\infty}\Big(\int_{\mathbb{R}^3}|f(s,y,\eta)|d\eta\Big)^{\f{p-1}{5p}}\right\}.
\end{align}
Then, using \eqref{5.24} and  similar arguments as in \eqref{5.26},  if $\v_0$ is small enough, one can obtain that 
\begin{align}\label{8.29}
&C_3\sup_{1\leq s\leq t,y\in\Omega}\left\{\|f(s)\|^{\f{4p+1}{5p}}_{L^\infty}\Big(\int_{\mathbb{R}^3}|f(s,y,\eta)|d\eta\Big)^{\f{p-1}{5p}}\right\}\leq \f12.
\end{align}
Substituting \eqref{8.29} into \eqref{8.25}, one proves that for $\b>\max\{\f92,4+|\g|\}$, 
\begin{align}\label{8.30}
\|f(t)\|_{L^\infty}\leq 2C_3\bar{M}^4(1+t)^{-r},\quad \forall\,t\geq0.
\end{align}
Taking 
$$\tilde{C}_3=2C_3\bar{M}^4,$$
then we obtain \eqref{1.30} from \eqref{8.30}.  Therefore the proof of Theorem \ref{thm1.3} is completed. $\hfill\Box$


\section{Appendix}

\subsection{Estimates on $K$}  In this subsection, we give the proof of some lemmas in Section 2 for completeness.

\medskip

\noindent{\bf Proof of Lemma \ref{lem2.1}:}  The estimate on $k_1(v,\eta)$ follows from a direct calculation.  We will mainly focus on $K_2(v,\eta)$.
It follows from \cite{Glassey,Bellomo} that 
\begin{align}\label{2.15-1}
0\leq k_2(v,\eta)=\f{c_2}{|v-\eta|^2}e^{-\f{|v-\eta|^2}{8}}e^{-\f{||v|^2-|\eta|^2|^2}{8|v-\eta|^2}}\int_{\mathbb{R}^2} B^{\ast}(|\eta-v|,|z_{\perp}|)e^{-\f{|z_{\perp}+\zeta_{\perp}|^2}{2}}dz_{\perp},
\end{align}
and $B^{\ast}(|\eta-v|,|z_{\perp}|)$ satisfies
\begin{equation}\label{2.16-1}
B^{\ast}(|\eta-v|,|z_{\perp}|)\leq C\f{|\eta-v|}{(|\eta-v|^2+|z_\perp|^2)^{\f{1-\g}{2}}}.
\end{equation}
where 
\begin{align}
&z=u-v,~~z_{\shortparallel}=[z\cdot\omega]\omega,~~z_{\perp}=z-z_{\shortparallel},~~\eta=v+z_{\shortparallel}.\notag
\end{align}
Then, substituting \eqref{2.16-1} into \eqref{2.15-1}, one obtains that 
\begin{align}
k_2(v,\eta)&\leq\f{c_2}{|v-\eta|}e^{-\f{|v-\eta|^2}{8}}e^{-\f{||v|^2-|\eta|^2|^2}{8|v-\eta|^2}}\int_{\mathbb{R}^2} (|\eta-v|^2+|z_{\perp}|^2)^{\f{\g-1}2}e^{-\f{|z_{\perp}+\zeta_{\perp}|^2}{2}}dz_{\perp}\nonumber\\
&\leq \f{c_2}{|v-\eta|^{\f{3-\g}{2}}}e^{-\f{|v-\eta|^2}{8}}e^{-\f{||v|^2-|\eta|^2|^2}{8|v-\eta|^2}}\int_{\mathbb{R}^2} |z_{\perp}|^{\f{\g-1}2}e^{-\f{|z_{\perp}+\zeta_{\perp}|^2}{2}}dz_{\perp}\nonumber\\
&\leq \f{C_\g}{|v-\eta|^{\f{3-\g}{2}}}e^{-\f{|v-\eta|^2}{8}}e^{-\f{||v|^2-|\eta|^2|^2}{8|v-\eta|^2}}.\notag
\end{align}
Thus  the proof of Lemma \ref{lem2.1} is completed. $\hfill\Box$

\

\noindent{\bf Proof Lemma \ref{lem2.2}:} Firstly, it is straightforward to prove \eqref{2.31} due to the fact that 
\begin{equation}
e^{-\f{|u|^2}{4}}\leq Ce^{-\f{|v|^2}{8}}~~~\mbox{for}~~|v-u|\leq 2m.\notag
\end{equation}
Next, we shall prove \eqref{2.33}. It is noted that
\begin{align}\label{2.35}
(k_1-k_1^m)(v,\eta)&=c|v-\eta|^\g [1-\chi_m(|v-u|)]e^{-\f{|v|^2}{4}}e^{-\f{|\eta|^2}{4}}\nonumber\\
&\leq c|v-\eta|^\g e^{-\f{|v|^2}{4}}e^{-\f{|\eta|^2}{4}}.
\end{align}
From \cite{Glassey}, we have that
\begin{align}\label{2.19}
k_2(v,\eta)&\leq \f{C}{|\eta-v|}e^{-\f{|\eta-v|^2}{8}}e^{-\f12|\zeta_{\shortparallel}|^2}\int_{\mathbb{R}^2}\f{1}{(|\eta-v|^2+|z_\perp|^2)^{\f{1-\g}{2}}}e^{-\f{|\zeta_{\perp}+z_{\perp}|^2}{2}}dz_{\perp},
\end{align}
where
\begin{align}\label{2.19-1}
&z=u-v,~~z_{\shortparallel}=[z\cdot\omega]\omega,~~z_{\perp}=z-z_{\shortparallel},~~\eta=v+z_{\shortparallel}\\ &\zeta:=
\f{1}{2}(v+\eta)=v+\f12z_{\shortparallel},~~\zeta_{\shortparallel}=[\zeta\cdot\omega]\omega.\label{2.19-2}
\end{align}
Denote
\begin{equation}\nonumber
\tilde{\chi}_m(s)=1-\chi_m(s)~~\mbox{for}~s\geq0.
\end{equation}
Then, it follows from \eqref{2.19} that  for $0\leq a\leq 1$, 
\begin{align}
&(k_2-k_2^m)(v,\eta)\nonumber\\
&\leq\f{1}{|\eta-v|}e^{-\f{|\eta-v|^2}{8}}e^{-\f12|\zeta_{\shortparallel}|^2}\int_{\mathbb{R}^2}\tilde\chi_m(\sqrt{|\eta-v|^2+|z_{\perp}|^2})(|\eta-v|^2+|z_{\perp}|^2)^{\f{\g-1}{2}}e^{-\f{|\zeta_{\perp}+z_{\perp}|^2}{2}}dz_{\perp}\nonumber\\
&\leq  \f{1}{|\eta-v|^{1+\f{(1-a)}{2}(1-\g)}}e^{-\f{|\eta-v|^2}{8}}e^{-\f12|\zeta_{\shortparallel}|^2}\int_{\mathbb{R}^2}\tilde\chi_m(\sqrt{|\eta-v|^2+|z_{\perp}|^2})\f{1}{(1+|\eta-v|^2+|z_{\perp}|^2)^{a\f{1-\g}{2}}}\nonumber\\
&\quad\times\f{(1+|\eta-v|^2+|z_{\perp}|^2)^{a\f{1-\g}{2}}}{(|\eta-v|^2+|z_{\perp}|^2)^{a\f{1-\g}{2}}}\cdot |z_{\perp}|^{(1-a)\f{\g-1}{2}} e^{-\f{|\zeta_{\perp}+z_{\perp}|^2}{2}}dz_{\perp}\nonumber\\
&\leq  C\f{1}{|\eta-v|^{1+\f{(1-a)}{2}(1-\g)}}e^{-\f{|\eta-v|^2}{8}}e^{-\f12|\zeta_{\shortparallel}|^2}\int_{\mathbb{R}^2}\tilde\chi_m(\sqrt{|\eta-v|^2+|z_{\perp}|^2})\f{1}{(1+|\eta-v|^2+|z_{\perp}|^2)^{a\f{1-\g}{2}}}\nonumber\\
&\quad\times\f{(1+|\eta-v|^2+|z_{\perp}|^2)^{a\f{1-\g}{2}}}{(m^2+|\eta-v|^2+|z_{\perp}|^2)^{a\f{1-\g}{2}}}\cdot |z_{\perp}|^{(1-a)\f{\g-1}{2}} e^{-\f{|\zeta_{\perp}+z_{\perp}|^2}{2}}dz_{\perp}\nonumber\\
&\leq  \f{m^{a(\g-1)}}{|\eta-v|^{1+\f{(1-a)}{2}(1-\g)}}e^{-\f{|\eta-v|^2}{8}}e^{-\f12|\zeta_{\shortparallel}|^2}\int_{\mathbb{R}^2}\f{|z_{\perp}|^{(1-a)\f{\g-1}{2}}}{(1+|\eta-v|^2+|z_{\perp}|^2)^{a\f{1-\g}{2}}} e^{-\f{|\zeta_{\perp}+z_{\perp}|^2}{2}}dz_{\perp}\nonumber\\
&\leq  \f{m^{a(\g-1)}}{|\eta-v|^{1+\f{(1-a)}{2}(1-\g)}}e^{-\f{|\eta-v|^2}{8}}e^{-\f12|\zeta_{\shortparallel}|^2}\int_{\mathbb{R}^2}\f{|\zeta_{\perp}-z_{\perp}|^{(1-a)\f{\g-1}{2}}}{(1+|\eta-v|^2+|\zeta_{\perp}-z_{\perp}|^2)^{a\f{1-\g}{2}}} e^{-\f{|z_{\perp}|^2}{2}}dz_{\perp},\notag
\end{align}
where in the last inequality we have made a change of variable $z_{\perp}+\zeta_{\perp}\rightarrow z_{\perp}$. Now we estimate  the RHS of above term. Following \cite{Strain-Guo}, we split it into two cases.

\medskip
\noindent{\it Case 1:} For $|z_{\perp}|\geq\f12|\zeta_\perp|$, it holds that 
\begin{align}
&(k_2-k_2^m)(v,\eta)\nonumber\\
&\leq  \f{m^{a(\g-1)}}{|\eta-v|^{1+\f{(1-a)}{2}(1-\g)}}e^{-\f{|\eta-v|^2}{8}}e^{-\f12|\zeta_{\shortparallel}|^2}\int_{\mathbb{R}^2}|\zeta_{\perp}-z_{\perp}|^{(1-a)\f{\g-1}{2}}e^{-\f{|z_{\perp}|^2}{4}}e^{-\f{|\z_{\perp}|^2}{16}}dz_{\perp}\nonumber\\
&\leq  \f{m^{a(\g-1)}}{|\eta-v|^{1+\f{(1-a)}{2}(1-\g)}}e^{-\f{|\eta-v|^2}{8}}e^{-\f14|\zeta_{\shortparallel}|^2}\int_{\mathbb{R}^2}|\zeta_{\perp}-z_{\perp}|^{(1-a)\f{\g-1}{2}}e^{-\f{|z_{\perp}|^2}{4}}e^{-\f{|\z|^2}{16}}dz_{\perp}\nonumber\\
&\leq  \f{m^{a(\g-1)}}{|\eta-v|^{1+\f{(1-a)}{2}(1-\g)}}e^{-\f{|\eta-v|^2}{10}}e^{-\f14|\zeta_{\shortparallel}|^2}\int_{\mathbb{R}^2}|\zeta_{\perp}-z_{\perp}|^{(1-a)\f{\g-1}{2}}e^{-\f{|z_{\perp}|^2}{4}}e^{-\f1{40}|\eta-v|^2-\f{|\z|^2}{16}}dz_{\perp}\nonumber.
\end{align}
It follows from \eqref{2.19-1} and \eqref{2.19-2} that
$$
|\eta-v|^2+|\zeta|^2=|\eta-v|^2+\f14|\eta+v|^2\geq \f14(|\eta|^2+|v|^2).
$$
This yields that
\begin{align}
&(k_2-k_2^m)(v,\eta)\nonumber\\
&\leq  \f{m^{a(\g-1)}}{|\eta-v|^{1+\f{(1-a)}{2}(1-\g)}}e^{-\f{|\eta-v|^2}{10}}e^{-\f14|\zeta_{\shortparallel}|^2}e^{-\f1{40}(|\eta|^2+|v|^2)}\int_{\mathbb{R}^2}|\zeta_{\perp}-z_{\perp}|^{(1-a)\f{\g-1}{2}}e^{-\f{|z_{\perp}|^2}{4}}dz_{\perp}\nonumber\\
&\leq  \f{C_\g m^{a(\g-1)}}{|\eta-v|^{1+\f{(1-a)}{2}(1-\g)}}e^{-\f{|\eta-v|^2}{10}}e^{-\f{||v|^2-|\eta|^2|^2}{16|v-\eta|^2}}e^{-\f1{40}(|\eta|^2+|v|^2)},\notag
\end{align}
where in the last inequality we have used the following fact
\begin{equation}\label{2.20}
|\zeta_{\shortparallel}|^2:=
|[\zeta\cdot\omega]\omega|^2= |[\zeta\cdot\omega]|^2= \Big|\zeta\cdot\f{z_{\shortparallel}}{|z_{\shortparallel}|}\Big|^2=\f14\Big|(\eta+v)\cdot\f{(\eta-v)}{|\eta-v|}\Big|^2=\f14\f{(|\eta|^2-|v|^2)^2}{|\eta-v|^2}.
\end{equation}

\noindent{\it Case 2:} For $|z_{\perp}|\leq\f12|\zeta_\perp|$, it holds that 
\begin{align}\label{2.39}
&(k_2-k_2^m)(v,\eta)\nonumber\\
&\leq  \f{C m^{a(\g-1)}}{|\eta-v|^{1+\f{(1-a)}{2}(1-\g)}}e^{-\f{|\eta-v|^2}{8}}e^{-\f12|\zeta_{\shortparallel}|^2}\int_{\mathbb{R}^2}\f{|\zeta_{\perp}-z_{\perp}|^{(1-a)\f{\g-1}{2}}}{(1+|\eta-v|^2+\f14|\zeta_{\perp}|^2)^{a\f{1-\g}{2}}} e^{-\f{|z_{\perp}|^2}{2}}dz_{\perp}\nonumber\\
&\leq  \f{C m^{a(\g-1)}}{|\eta-v|^{1+\f{(1-a)}{2}(1-\g)}}e^{-\f{|\eta-v|^2}{8}}e^{-\f14|\zeta_{\shortparallel}|^2}\int_{\mathbb{R}^2}\f{|\zeta_{\perp}-z_{\perp}|^{(1-a)\f{\g-1}{2}}}{(1+|\eta-v|^2+|\zeta_{\shortparallel}|^2+|\zeta_{\perp}|^2)^{a\f{1-\g}{2}}} e^{-\f{|z_{\perp}|^2}{2}}dz_{\perp}\nonumber\\
&\leq  \f{C m^{a(\g-1)}}{|\eta-v|^{1+\f{(1-a)}{2}(1-\g)}}e^{-\f{|\eta-v|^2}{8}}e^{-\f14|\zeta_{\shortparallel}|^2}\int_{\mathbb{R}^2}\f{|\zeta_{\perp}-z_{\perp}|^{(1-a)\f{\g-1}{2}}}{(1+|\eta-v|^2+|\zeta|^2)^{a\f{1-\g}{2}}} e^{-\f{|z_{\perp}|^2}{2}}dz_{\perp}\nonumber\\
&\leq  \f{C m^{a(\g-1)}}{|\eta-v|^{1+\f{(1-a)}{2}(1-\g)}}\f{1}{(1+|v|^2+|\eta|^2)^{a\f{1-\g}{2}}}e^{-\f{|\eta-v|^2}{8}}e^{-\f14|\zeta_{\shortparallel}|^2}\int_{\mathbb{R}^2}|\zeta_{\perp}-z_{\perp}|^{(1-a)\f{\g-1}{2}} e^{-\f{|z_{\perp}|^2}{2}}dz_{\perp}\nonumber\\
&\leq  \f{C_\g m^{a(\g-1)}}{|\eta-v|^{1+\f{(1-a)}{2}(1-\g)}}\f{1}{(1+|v|+|\eta|)^{a(1-\g)}}e^{-\f{|\eta-v|^2}{8}}e^{-\f{||v|^2-|\eta|^2|^2}{16|v-\eta|^2}},
\end{align}
where we have used \eqref{2.20} in the last inequality. Combining \eqref{2.35}-\eqref{2.39}, one proves \eqref{2.33}. Therefore,  the proof of Lemma \ref{lem2.2} is completed. $\hfill\Box$

\

\subsection{Local-in-time existence}

In the following, we consider the local existence of unique solutions to the Boltzmann equation \eqref{1.1} with large initial data  in $L^\infty$-norm.

\medskip

\noindent{\bf Proof of Proposition \ref{prop7.1}:} To prove the local existence for the Boltzmann equation, we consider the iteration that for $n=0,1,2,\cdots$, 
\begin{align}\label{7.2}
& F^{n+1}_t+v\cdot\nabla_x F^{n+1}+F^{n+1}\cdot\int_{\mathbb{R}^3}\int_{\mathbb{S}^2}B(v-u,\t)F^{n}(t,x,u)d\omega du=Q_+(F^n,F^n)
\end{align}
with 
\begin{align}\label{7.3}
F^{n+1}(t,x,v)\Big|_{t=0}=F_0(x,v)\geq0,~~~~\mbox{and}~~~F^0(t,x,v)=\mu(v).
\end{align}
Denote 
\begin{equation}
f^{n+1}=\f{F^{n+1}-\mu}{\sqrt{\mu}}.\notag
\end{equation}
Then \eqref{7.2} can be written equivalently as 
\begin{align}\label{7.5}
&f^{n+1}_t+v\cdot\nabla_x f^{n+1}+f^{n+1}\cdot \int_{\mathbb{R}^3}\int_{\mathbb{S}^2}B(v-u,\t)\Big\{\mu(u)+\sqrt{\mu}f^n(t,x,u)\Big\}d\omega du\nonumber\\
&=Kf^n+\f{1}{\sqrt{\mu}}Q_+(\sqrt{\mu}f^n,\sqrt{\mu}f^n)
\end{align}
with $n=0,1,2,\cdots$ and 
\begin{align}\label{7.6}
f^{n+1}(t,x,v)\Big|_{t=0}=f_0(x,v)~~~\mbox{and}~~~f^0(t,x,v)=0.
\end{align}
It is a normal procedure to  solve the approximated problems \eqref{7.2}-\eqref{7.3} (or equivalently \eqref{7.5}-\eqref{7.6}) since they are linear at each step and the angular cutoff assumption is posed. Then we get an approximation sequence $F^{n+1},~n=0,1,2,\cdots$.

Firstly, we consider the positivity of $F^{n+1}$.  It is noted that 
\begin{align}\label{7.7}
F^{n+1}(t,x,v)&=e^{-\int_0^tg^n(\tau,x-v(t-\tau),v)d\tau}\cdot F_0(x-vt,v)\nonumber\\
&\quad+\int_0^te^{-\int_s^tg^n(\tau,x-v(t-\tau),v)d\tau}\cdot Q_+(F^n,F^n)(s,x-v(t-s),v)ds,
\end{align}
with
\begin{align}
\notag
g^n(\tau,y,v)&=\int_{\mathbb{R}^3}\int_{\mathbb{S}^2}B(v-u,\theta)F^n(\tau,y,u)dud\omega\nonumber\\
&=\int_{\mathbb{R}^3}\int_{\mathbb{S}^2}B(v-u,\t)\Big\{\mu(u)+\sqrt{\mu}f^n(\tau,y,u)\Big\}d\omega du.\notag
\end{align}
By induction on $n$, we can prove that if $F^n\geq0$, then it holds that  
$$g^n(s,x-v(t-s),v)\geq0~~\mbox{and}~~Q_+(F^n,F^n)(s,x-v(t-s),v)\geq0,$$
which together with \eqref{7.7} yield that 
\begin{align}
\notag
F^{n+1}(t,x,v)&\geq e^{-\int_0^tg^n(\tau,x-v(t-\tau),v)d\tau}\cdot F_0(x-vt,v)\geq 0.
\end{align}
Therefore, we have proved the positivity of the approximation sequences, i.e., $F^{n+1}\geq0,~n=0,1,\cdots$.

Next, we consider the  uniform estimate for the approximation sequence.  And it is more convenient to use the equivalent form $f^{n+1}$. Then, it follows from \eqref{7.5} that 
\begin{align}\label{7.10}
& f^{n+1}(t,x,v)\nonumber\\
&=e^{-\int_0^tg^n(\tau,x-v(t-\tau),v)d\tau}\cdot f_0(x-vt,v) \nonumber\\
&\quad+\int_0^te^{-\int_s^tg^n(\tau,x-v(t-\tau),v)d\tau}\cdot (Kf^n)(s,x-v(t-s),v)ds\nonumber\\
&\quad+\int_0^te^{-\int_s^tg^n(\tau,x-v(t-\tau),v)d\tau}\cdot\f{1}{\sqrt{\mu(v)}}Q_+(\sqrt{\mu}f^n,\sqrt{\mu}f^n)(s,x-v(t-s),v)ds,
\end{align}
which yields that 
\begin{align}\label{7.11}
|w_\b(v) f^{n+1}(t,x,v)|&\leq \|w_\b(v)f_0\|_{L^\infty} +\int_0^t\Big|w_\b(v) (Kf^n)(s,x-v(t-s),v)\Big| ds\nonumber\\
&\quad+\int_0^t\f{w_\b(v)}{\sqrt{\mu(v)}}\Big|Q_+(\sqrt{\mu}f^n,\sqrt{\mu}f^n)(s,x-v(t-s),v)\Big|ds.
\end{align}
It follows from \eqref{2.18} that 
\begin{align}\label{7.13}
\int_0^t\Big|w_\b(v) (Kf^n)(s,x-v(t-s),v)\Big|ds
&\leq \int_0^t\|w_\b f^n(s)\|_{L^\infty}ds\int_{\mathbb{R}^3}w_\b(v)|k(v,\eta)|w_\b (\eta)^{-1}d\eta\nonumber\\
&\leq  C_\g\int_0^t\|w_\b f^n(s)\|_{L^\infty}ds.
\end{align}
To estimate the last term on the RHS of \eqref{7.11}, by similar arguments as in \eqref{4.15}, one  gets that 
\begin{align}\label{7.15}
&\f{w_\b (v)}{\sqrt{\mu(v)}}\Big|Q_+(\sqrt{\mu}f^n,\sqrt{\mu}f^n)(s,x-v(t-s),v)\Big|\nonumber\\
&\leq   C\int_{\mathbb{R}^3}\int_{\mathbb{S}^2}B(v-u,\t)\sqrt{\mu(u)}\Big|w_\b(u')f^n(s,x-v(t-s),u')f^n(s,x-v(t-s),v')
\Big| dud\omega\nonumber\\
&\quad+C\int_{\mathbb{R}^3}\int_{\mathbb{S}^2}B(v-u,\t)\sqrt{\mu(u)}\Big|f^n(s,x-v(t-s),u')w_\b(v')f^n(s,x-v(t-s),v')
\Big|dud\omega,\nonumber\\
&:=
I_1+I_2.
\end{align}
It follows from the change of variables \eqref{3.9}-\eqref{3.11} that for $\b>3$,  
\begin{align}\label{7.20}
I_{1}
&\leq C\|w_\b f^n(s)\|_{L^\infty}\int_{\mathbb{R}^3}\int_{\mathbb{S}^2}\f{|z_{\shortparallel}|}{(|z_{\shortparallel}|+|z_{\perp}|)^{1-\g}}
e^{-\f{|v+z|^2}{4}}\Big|f^n(s,x-v(t-s),v+z_{\shortparallel}) \Big| dz_{\perp}d|z_{\shortparallel}|d\omega\nonumber\\
&\leq C\|w_\b f^n(s)\|_{L^\infty}\int_{\mathbb{R}^3}\int_{z_{\perp}}\f{|z_{\perp}|^{\f{\g-1}{2}}}{|\eta-v|^{\f{3-\g}{2}}}
e^{-\f{|\eta+z_{\perp}|^2}{4}}\Big|f^n(s,x-v(t-s),\eta) \Big| dz_{\perp}d\eta\nonumber\\
&\leq C_\g\|w_\b f^n(s)\|_{L^\infty} \int_{\mathbb{R}^3}\f{1}{|\eta-v|^{\f{3-\g}{2}}}
\Big|f^n(s,x-v(t-s),\eta) \Big| d\eta \nonumber\\
&\leq C_\g\|w_\b f^n(s)\|^2_{L^\infty} \int_{\mathbb{R}^3}\f{(1+|\eta|)^{-\b}}{|\eta-v|^{\f{3-\g}{2}}}
d\eta \leq C(\g,\b)\|w_\b f^n(s)\|^2_{L^\infty}.
\end{align}
It is noted that by a rotation, one obtains the interchange of $v'$ and $u'$. Then, one can change $I_{2}$  to a similar form as $I_{1}$.  Thus, by  similar arguments as above, one can obtain that  for $\b>3$,
\begin{align}
\notag
I_{2} \leq C(\g,\b)\|w_\b f^n(s)\|^2_{L^\infty},
\end{align}
which together with \eqref{7.20} and \eqref{7.15} yields that  for $\b>3$, 
\begin{equation}\label{7.22}
\int_0^t\f{w_\b (v)}{\sqrt{\mu(v)}}\Big|Q_+(\sqrt{\mu}f^n,\sqrt{\mu}f^n)(s,x-v(t-s),v)\Big|ds
\leq C(\g,\b)\int_0^t\|w_\b f^n(s)\|^2_{L^\infty}ds.
\end{equation}
Substituting \eqref{7.22} and \eqref{7.13} into \eqref{7.11}, one obtains that  for $\b>3$,
\begin{equation}\label{7.23}
\|w_{\b}(v)f^{n+1}(t)\|_{L^\infty}\leq \|w_{\b}f_0\|_{L^\infty}+C_4t\Big\{\sup_{0\leq s\leq t}\|w_{\b}f^n(s)\|_{L^\infty}+\sup_{0\leq s\leq t}\|w_{\b}f^n(s)\|^2_{L^\infty}\Big\},
\end{equation}
where the positive constant $C_4\geq1$ depends only on $\g,\b$.  By induction on $n$, we can prove that if 
\begin{align}\label{7.24}
\sup_{0\leq s\leq t_1}\|w_{\b}f^n(s)\|_{L^\infty}\leq 2\|w_{\b}f_0\|_{L^\infty},\quad t_1=\Big(8C_4[1+\|w_{\b}f_0\|_{L^\infty}]\Big)^{-1},
\end{align}
then it follows from \eqref{7.23}  and \eqref{7.24} that for $\b>3$,
\begin{equation}\label{7.25}
\sup_{0\leq s\leq t_1}\|w_{\b}f^{n+1}(s)\|_{L^\infty}\leq 2\|w_{\b}f_0\|_{L^\infty}~\mbox{with}~t_1=\Big(8C_4[1+\|w_{\b}f_0\|_{L^\infty}]\Big)^{-1},
\end{equation}
for all $n\geq 0$.

Now we  prove that $f^{n+1},~n=0,1,2,\cdots$ is a Cauchy sequence.  It follows from \eqref{7.10} that 
\begin{align}\label{7.26}
&\Big|\sqrt{w_{\b}(v)}(f^{n+2}-f^{n+1})(t,x,v)]\Big|\nonumber\\
&\leq |\sqrt{w_{\b}(v)}f_0(x-vt,v)|\cdot\int_0^t\Big|(g^{n+1}-g^n)(\tau,x-v(t-\tau),v)\Big|d\tau\nonumber\\
&\quad+\int_0^t|\sqrt{w_{\b}(v)}Kf^{n+1}(s,x-v(t-s),v)|\cdot\int_s^t\Big|(g^{n+1}-g^n)(\tau,x-v(t-\tau),v)\Big|d\tau ds\nonumber\\
&\quad+\int_0^t\f{\sqrt{w_{\b}(v)}}{\sqrt{\mu(v)}}|Q_+(\sqrt{\mu}f^{n+1},\sqrt{\mu}f^{n+1})(s,x-v(t-s),v)|\nonumber\\
&\qquad\qquad\qquad\times\int_s^t\Big|(g^{n+1}-g^n)(\tau,x-v(t-\tau),v)\Big|d\tau ds\nonumber\\
&\quad+\int_0^t|\sqrt{w_{\b}(v)}(Kf^{n+1}-Kf^n)(s,x-v(t-s),v)|ds\nonumber\\
&\quad+\int_0^t\f{\sqrt{w_{\b}(v)}}{\sqrt{\mu(v)}}\Big|\Big(Q_+(\sqrt{\mu}f^{n+1},\sqrt{\mu}f^{n+1})-Q_+(\sqrt{\mu}f^{n},\sqrt{\mu}f^{n})\Big)(s,x-v(t-s),v)\Big|ds\nonumber\\
&:=
I_3+I_4+I_5+I_6+I_7.
\end{align}
A direct calculation shows that 
\begin{align}
\Big|(g^{n+1}-g^n)(\tau,y,u)\Big|
\leq C\nu(v)\|(f^{n+1}-f^{n})(\tau)\|_{L^\infty},\notag
\end{align}
which, together with \eqref{7.13} and \eqref{7.22}, yields that, for $0\leq t\leq t_1$,
\begin{align}\label{7.27}
I_3+I_4+I_5&\leq C\Big\{ \|\nu\sqrt{w_{\b}} f_0\|_{L^\infty}+\int_0^t\Big|\sqrt{w_{\b}(v)}\nu(v)Kf^{n+1}(s,x-v(t-s),v)\Big|ds\nonumber\\
&\quad+\int_0^t\f{\sqrt{w_{\b}(v)}\nu(v)}{\sqrt{\mu(v)}}\Big|Q_+(\sqrt{\mu}f^{n+1},\sqrt{\mu}f^{n+1})(s,x-v(t-s),v)\Big|ds\Big\}\nonumber\\
&\qquad\qquad\times\int_0^t\|(f^{n+1}-f^{n})(\tau)\|_{L^\infty}d\tau\nonumber\\
&\leq Ct\Big\{  \|w_{\b}f_0\|_{L^\infty}+t\sup_{0\leq s\leq t}\|w_{\b}f^{n+1}(s)\|_{L^\infty}+t\sup_{0\leq s\leq t}\|w_{\b}f^{n+1}(s)\|^2_{L^\infty}\Big\}\nonumber\\
&\qquad\qquad\times\sup_{0\leq \tau\leq t}\|(f^{n+1}-f^{n})(\tau)\|_{L^\infty}\nonumber\\
&\leq Ct\|w_{\b}f_0\|_{L^\infty}\cdot\sup_{0\leq \tau\leq t}\|(f^{n+1}-f^{n})(\tau)\|_{L^\infty},
\end{align}
with $\b>3$, where we have used the uniform estimate \eqref{7.25} in the last inequality.
By the same arguments as \eqref{7.13}, one can obtains that 
\begin{align}\label{7.28}
I_6\leq Ct\sup_{0\leq \tau\leq t}\|\sqrt{w_{\b}}(f^{n+1}-f^{n})(\tau)\|_{L^\infty}.
\end{align}
To prove $I_7$, we note that 
\begin{align}\label{7.29}
&\Big|\Big(Q_+(\sqrt{\mu}f^{n+1},\sqrt{\mu}f^{n+1})-Q_+(\sqrt{\mu}f^{n},\sqrt{\mu}f^{n})\Big)(s,x-v(t-s),v)\Big|\nonumber\\
&\leq \Big|Q_+(\sqrt{\mu}f^{n+1},\sqrt{\mu}(f^{n+1}-f^n))(s,x-v(t-s),v)\Big|\nonumber\\
&\quad+\Big|Q_+(\sqrt{\mu}(f^{n+1}-f^n),\sqrt{\mu}f^{n})(s,x-v(t-s),v)\Big|.
\end{align}
Denoting $y=x-v(t-s)$, by similar arguments as in \eqref{7.15},  we have that 
\begin{align}\label{7.30}
&\f{\sqrt{w_{\b}(v)}}{\sqrt{\mu(v)}}\Big|Q_+(\sqrt{\mu}f^{n+1},\sqrt{\mu}(f^{n+1}-f^n))(s,x-v(t-s),v)\Big|\nonumber\\
&\leq C\int_{\mathbb{R}^3}\int_{\mathbb{S}^2}B(v-u,\t)\sqrt{\mu(u)}|f^{n+1}(s,y,u')|\cdot\Big|\sqrt{w_{\b}(v')}(f^{n+1}-f^n)(s,y,v')\Big|dud\omega\nonumber\\
&\quad+C\f{\nu(v)}{\sqrt{w_\b(v)}}\|w_{\b}f^{n+1}(s)\|_{L^\infty}\cdot\|\sqrt{w_{\b}}(f^{n+1}-f^n)(s)\|_{L^\infty}.
\end{align}
To estimate the first term on the RHS of \eqref{7.30}, by a rotation, we  interchange $v'$ and $u'$, then using the same arguments as in \eqref{3.9}-\eqref{7.20}, one can obtain that 
\begin{align}
& C\int_{\mathbb{R}^3}\int_{\mathbb{S}^2}B(v-u,\t)\sqrt{\mu(u)}|f^{n+1}(s,y,u')|\cdot\Big|\sqrt{w_{\b}(v')}(f^{n+1}-f^n)(s,y,v')\Big|dud\omega\nonumber\\
&\leq C\|\sqrt{w_{\b}}(f^{n+1}-f^n)(s)\|_{L^\infty}\int_{\mathbb{R}^3}\int_{z_{\perp}}\f{|z_{\perp}|^{\f{\g-1}{2}}}{|\eta-v|^{\f{3-\g}{2}}}
e^{-\f{|\eta+z_{\perp}|^2}{4}}\Big|f^{n+1}(s,y,\eta) \Big| dz_{\perp}d\eta\nonumber\\
&\leq C\|w_{\b}f^{n+1}(s)\|_{L^\infty} \cdot\|\sqrt{w_{\b}}(f^{n+1}-f^n)(s)\|_{L^\infty},\notag
\end{align}
with $\b>3$, which together with \eqref{7.30} yield that 
\begin{align}\label{7.32}
&\f{\sqrt{w_{\b}(v)}}{\sqrt{\mu(v)}}\Big|Q_+(\sqrt{\mu}f^{n+1},\sqrt{\mu}(f^{n+1}-f^n))(s,x-v(t-s),v)\Big|\nonumber\\
&\leq C\|w_{\b}f^{n+1}(s)\|_{L^\infty}\cdot\|\sqrt{w_{\b}}(f^{n+1}-f^n)(s)\|_{L^\infty},
\end{align}
for $\b>3$. 
Similarly, one can gets that for $\b>3$,
\begin{align}\label{7.33}
&\f{\sqrt{w_{\b}(v)}}{\sqrt{\mu(v)}}\Big|Q_+(\sqrt{\mu}(f^{n+1}-f^n),\sqrt{\mu}f^{n})(s,x-v(t-s),v)\Big|\nonumber\\
&\leq C\|w_{\b}f^{n}(s)\|_{L^\infty}\cdot\|\sqrt{w_{\b}}(f^{n+1}-f^n)(s)\|_{L^\infty}.
\end{align}
Then it follows from \eqref{7.29}, \eqref{7.32} and \eqref{7.33} that for $0\leq t\leq t_1$,
\begin{align}\label{7.34}
I_7&\leq Ct\sup_{0\leq s\leq t}\Big\{\|w_{\b}f^{n+1}(s)\|_{L^\infty}+\|w_{\b}f^{n}(s)\|_{L^\infty}\Big\}\cdot\sup_{0\leq s\leq t}\|\sqrt{w_{\b}}(f^{n+1}-f^n)(s)\|_{L^\infty}\nonumber\\
&\leq Ct\|w_{\b}f_0\|_{L^\infty}\sup_{0\leq s\leq t}\|\sqrt{w_{\b}}(f^{n+1}-f^n)(s)\|_{L^\infty},
\end{align}
with $\b>3$, 
where we have used the uniform estimate \eqref{7.25} in the last inequality.

Substituting  \eqref{7.27}, \eqref{7.28} and  \eqref{7.34} into \eqref{7.26}, one obtains that for $0\leq t\leq t_1$, 
\begin{align}\notag
\sup_{0\leq s\leq T_1}\|\sqrt{w_{\b}}(f^{n+2}-f^{m+1})(s)\|&\leq Ct_1(1+\|w_{\b}f_0\|_{L^\infty})\cdot\sup_{0\leq s\leq T_1}\|\sqrt{w_{\b}}(f^{n+1}-f^n)(s)\|_{L^\infty}\nonumber\\
&\leq \f{C}{8C_4}\sup_{0\leq s\leq t_1}\|\sqrt{w_{\b}}(f^{n+1}-f^n)(s)\|_{L^\infty}\nonumber\\
&\leq \f12\sup_{0\leq s\leq t_1}\|\sqrt{w_{\b}}(f^{n+1}-f^n)(s)\|_{L^\infty},\notag
\end{align}
where we have chosen $C_4$ suitably large such that $\f{C}{8C_4}\leq \f12$. Thus, by induction on $n$, it is direct to obtain that 
\begin{align}\notag
\sup_{0\leq s\leq t_1}\|\sqrt{w_{\b}}(f^{n+2}-f^{n+1})(s)\|\leq 2^{-n-1}\|\sqrt{w_{\b}}(f^1-f^0)\|_{L^\infty}
\leq 2^{-n}\|w_{\b}f_0\|_{L^\infty},
\end{align}
which yields immediately that $f^{n+1},~n=0,1,2,\cdots$ is a Cauchy sequence. Therefore, there exists a limit $f$ such that 
\begin{equation}\notag
\sup_{0\leq s\leq t_1}\|\sqrt{w_{\b}}(f^n-f)(s)\|_{L^\infty}\rightarrow0 ~~\mbox{as}~~n\rightarrow +\infty.
\end{equation}
The limit function $f$ is indeed a mild solution to the Boltzmann equation \eqref{1.1} and \eqref{1.5-1}. It follows from \eqref{7.25} that 
\begin{equation}\notag
\sup_{0\leq t\leq t_1}\|w_{\b}f(t)\|_{L^\infty}\leq 2\|w_{\b}f_0\|_{L^\infty}.
\end{equation}

Now we consider the uniqueness. Let $\tilde{f}(t,x,v)$ be another solution of the Boltzmann equation \eqref{1.1} and \eqref{1.5-1} with the bound $\sup_{0\leq t\leq t_1}\|w_{\b}\tilde{f}(t)\|_{L^\infty}<+\infty$, by similar arguments as \eqref{7.26}-\eqref{7.34},  it is direct to obtain that 
\begin{align}\notag
\|\sqrt{w_{\b}}(f-\tilde{f})(t)\|_{L^\infty}\leq C(1+\|w_{\b}f\|_{L^\infty}+\|w_{\b}\tilde{f}\|_{L^\infty})\cdot\int_0^t\|\sqrt{w_{\b}}(f-\tilde{f})(s)\|_{L^\infty}ds,
\end{align}
which together with the Gronwall inequality yields the uniqueness, i.e., $f=\tilde{f}$.

Multiplying  \eqref{7.2}  by  $1,v,|v|^2$ and $F^{m+1}$,  integrating by parts and then taking the limit  $m\rightarrow+\infty$, one can directly obtain \eqref{1.13}-\eqref{1.14-1}.  

Finally, if $F_0$ (or equivalent $f_0$) is continuous, it is direct to check that $F^{n+1}(t,x,v)$ (or equivalent $f^{n+1}(t,x,v)$) is continuous in $[0,\infty)\times\Omega\times\mathbb{R}^3$. The continuous of $f(t,x,v)$  is an immediate consequence of  $\sup_{0\leq s\leq t_1}\|(f^{n+1}-f)(s)\|_{L^\infty}\rightarrow 0$ as $n\rightarrow +\infty$. Therefore the proof of  Proposition \ref{prop7.1} is completed. $\hfill\Box$

\bigskip

\noindent {\bf Acknowledgments:} Renjun Duan is partially supported by the General Research Fund (Project No. 409913) from RGC of Hong Kong. Feimin Huang is partially supported by National Center for Mathematics
and Interdisciplinary Sciences, AMSS, CAS and NSFC Grant No. 11371349.  Yong Wang is partially supported by National Natural
Sciences Foundation of China No.  11401565. Tong Yang is partially supported by the
General Research Fund of Hong Kong, CityU 103412.



\begin{thebibliography}{99}
	
	


	

\bibitem{BM} 
C. Baranger, C. Mouhot, 
{\it Explicit spectral gap estimates for the linearized Boltzmann and Landau operators with hard potentials},
Rev. Mat. Iberoamericana {\bf 21} (2005), no. 3, 819-841.

\bibitem{Bellomo} 
\newblock N. Bellomo, A. Palczewski, G. Toscani,
\emph{ Mathematical topics in nonlinear kinetic theory,}
\newblock  World Scientific Publishing Co., Singapore, 1988. 



\bibitem{BG}
M. Briant, Y. Guo,
{\it Asymptotic stability of the Boltzmann equation with Maxwell boundary conditions},
arXiv:1511.01305.




%

\bibitem{Car}
T. Carleman, 
{\it Sur la th\'eorie de l'\'equation int\'egrodiff\'erentielle de Boltzmann}, 
Acta Math. {\bf 60} (1933), no. 1, 91-146.



\bibitem{CIP}
\newblock C. Cercignani, R. Illner, M. Pulvirenti, 
\emph{The Mathematical Theory of Dilute Gases},  Springer-Verlag, New York, 1994.


\bibitem{Duan-Yang-Zhao}
\newblock R. J. Duan, T. Yang, H. J. Zhao, 
\emph{ The Vlasov-Poisson-Boltzmann system for soft potentials, }
\newblock Mathematical Models and Methods in Applied Sciences.  \textbf{23}(2013), no. 6, 979-1028.

\bibitem{D-Lion} 
\newblock R.J. DiPerna, P.-L. Lions, 
\emph{On the Cauchy problem for Boltzmann equation: Global existence and weak stability, }
\newblock  Ann. of Math. \textbf{130} (1989), 321-366.

\bibitem{Desvillettes-V} 
\newblock L. Desvillettes,  C. Villani, 
\emph{On the trend to global equilibrium for spatially inhomogeneous kinetic systems:The Boltzmann equation,}
\newblock Invent. Math. \textbf{159} (2005), 243-316.


\bibitem{EP}
R.  Ellis, M.A. Pinsky, 
{\it The first and second fluid approximations to the linearized Boltzmann equation},
 J. Math. Pures Appl. (9) {\bf 54} (1975), 125-156.





\bibitem{HW}
\newblock F. M. Huang, Y. Wang,
\emph{Macroscopic regularity for the Boltzmann equation,}
\newblock Preprinted  2015.

\bibitem{Glassey} 
\newblock R.T. Glassey,
\emph{The cauchy problem in kinetic theory,}
\newblock Society for Industrial and Applied Mathematics (SIAM), Philadelphia,
1996.

\bibitem{Gra}
H. Grad,
{\it Asymptotic theory of the Boltzmann equation}, in Rarefied Gas Dynamics, edited by J. A. Laurmann,  (Academic Press, New York, 1963), Vol. 1, pp. 26-59.


\bibitem{GMM} 
M. P. Gualdani, S. Mischler, C. Mouhot,
{\it Factorization for non-symmetric operators and exponential H-theorem},
arXiv:1006.5523.

\bibitem{Guo-03}
Y. Guo, 
{\it Classical solutions to the Boltzmann equation for molecules with an angular cutoff},
 Arch. Ration. Mech. Anal. {\bf 169} (2003), no. 4, 305-353.
%



\bibitem{Guo2}
\newblock Y. Guo,
\emph{Decay and continuity of the Boltzmann equation in Bounded domains},
\newblock  Arch. Rational. Mech. Anal. \textbf{197}(2010), 713-809.

\bibitem{Guo}
\newblock Y. Guo,
\emph{Bounded solutions for the Boltzmann equation},
\newblock  Quart. Appl. Math. \textbf{68} (2010), no.1, 143-148.


\bibitem{IIIner-S} 
\newblock R. Illner, M. Shinbrot, 
\emph{Global existence for a rare gas in an infinite vacuum, }
\newblock Comm. Math. Phys. \textbf{95} (1984), 117-126.


\bibitem{Kaniel-S} 
\newblock S. Kaniel, M. Shinbrot, 
\emph{The Boltzmann equation I: Uniqueness and local existence, }
\newblock Comm. Math. Phys. \textbf{58} (1978) 65-84.

\bibitem{Kim2014}
\newblock C. Kim,
\emph{Boltzmann equation with a large potential in a periodic box},
\newblock Comm. Partial Diff. Eqns.  \textbf{39}(2014), 1393-1423.


\bibitem{Liu-Yang-Yu}
\newblock T. Liu, T. Yang,  S. H. Yu,
 \emph{Energy method for the Boltzmann equation,}
\newblock  Physica D, \textbf{188}(2004), 178-192.

\bibitem{LM15}
X.-G. Lu, C. Mouhot, 
{\it On measure solutions of the Boltzmann equation, Part II: Rate of convergence to equilibrium}, J. Differential Equations {\bf 258} (2015), no. 11, 3742-3810.

%






%

\bibitem{Strain}
\newblock R.M. Strain, 
\emph{Optimal time decay of the non cut-off Boltzmann equation in the whole space},
\newblock  Kinetic and Related Models. \textbf{5}(2012), 583-613.

\bibitem{Strain-Guo}
\newblock R.M. Strain,  Y. Guo,
\emph{Exponential decay for soft potentials near Maxwellian},
\newblock  Arch. Rational. Mech. Anal. \textbf{187}(2008), 287-339.

\bibitem{Uk}
S. Ukai,
{\it On the existence of global solutions of mixed problem for non-linear Boltzmann equation}, Proc. Jpn. Acad. {\bf 50} (1974), 179-184.

\bibitem{Ukai-Yang}
\newblock S. Ukai and  T. Yang,
\emph{The Boltzmann equation in the space $L^2\cap L^\infty_\beta$ : Global and time-periodic solutions,} 
\newblock Anal. Appl. \textbf{4} (2006),  263-310.

\bibitem{Vidav}
\newblock I. Vidav,
\emph{Spectra of perturbed semigroups with applications to transport theory,}
\newblock J. Math. Anal. Appl., \textbf{30}(1970), 264-279.

\bibitem{Vi} 
\newblock C. Villani,  
\emph{A review of mathematical topics in collisional kinetic theory. Handbook of mathematical
fluid dynamics, }
\newblock Vol. I, 71-305. North-Holland, Amsterdam, 2002.


\end{thebibliography}
\end{document}